\documentclass[a4paper,12pt,reqno,oneside]{amsart} 

\usepackage{setspace} 
\usepackage{typearea} 
\usepackage{amscd,amsmath,amssymb,verbatim} 

\usepackage[dvipsnames]{xcolor}
\usepackage{xr-hyper}
\usepackage{cite}
\usepackage[colorlinks,backref,final,bookmarksnumbered,bookmarks]{hyperref}
\usepackage[backrefs]{amsrefs}
\hypersetup{linkcolor=bluex,citecolor=Green,filecolor=cyan,urlcolor=magenta}
\usepackage{graphicx}
\definecolor{bluex}{Hsb}{220, 1, 1}

\newcommand{\borderref}[2]{{\hypersetup{linkcolor=black}\hyperref[#1]{#1 #2}}}
  
\usepackage{ifthen}
\newcommand{\arxiv}[2][]{\ifthenelse{\equal{#1}{}}
{\href{http://arxiv.org/abs/#2}{\tt arXiv:#2}}
{\href{http://arxiv.org/abs/math/#2}{\tt arXiv:math.#1/#2}}}

\theoremstyle{plain}
\newtheorem{theorem}{Theorem}[section]
\newtheorem{mainthm}{Theorem}
 
\newtheorem{lemma}[theorem]{Lemma}
\newtheorem{proposition}[theorem]{Proposition}
\newtheorem{corollary}[theorem]{Corollary}
\newtheorem*{corollary*}{Corollary}

\theoremstyle{definition}
\newtheorem{remark}[theorem]{Remark}
\newtheorem*{remark*}{Remark}
\newtheorem{example}[theorem]{Example}

\def\x{\times}
\def\but{\setminus}

\def\eps{\varepsilon}
\def\phi{\varphi}
\def\emptyset{\varnothing}
\renewcommand{\:}{\colon}

\def\R{\mathbb{R}}
\def\Q{\mathbb{Q}}
\def\Z{\mathbb{Z}}

\def\L{\mathcal{L}}
\def\xr#1{\xrightarrow{#1}} 

\makeatletter
\newcommand{\xR}[2][]{\ext@arrow 0359\Rightarrowfill@{#1}{#2}}
\newcommand{\xL}[2][]{\ext@arrow 0359\Leftarrowfill@{#1}{#2}}
\makeatother
\DeclareMathOperator{\lk}{lk}

\DeclareMathOperator{\rk}{rk}
\DeclareMathOperator{\Hom}{Hom}
\DeclareMathOperator{\adj}{adj}
\DeclareMathOperator{\sgn}{sgn}

\newcommand{\newsym}[5]{\fontfamily{#2}\fontencoding{#1}\fontseries{#3}\fontshape{#4}\selectfont\char#5}
\makeatletter
\newcommand{\newmathsymbol}[6]{#1{\@Pimathsymbol{#2}{#3}{#4}{#5}{#6}}}
\def\@Pimathsymbol#1#2#3#4#5{\mathchoice
  {\@Pim@thsymbol{#1}{#2}{#3}{#4}{#5}\tf@size}
  {\@Pim@thsymbol{#1}{#2}{#3}{#4}{#5}\tf@size}
  {\@Pim@thsymbol{#1}{#2}{#3}{#4}{#5}\sf@size}
  {\@Pim@thsymbol{#1}{#2}{#3}{#4}{#5}\ssf@size}}
\def\@Pim@thsymbol#1#2#3#4#5#6{\mbox{\fontsize{#6}{#6}\newsym{#1}{#2}{#3}{#4}{#5}}}
\makeatother

\DeclareFontEncoding{LS1}{}{}
\DeclareFontSubstitution{LS1}{stix}{m}{n}
\def\varnabla{\newmathsymbol{\mathord}{LS1}{stix}{m}{n}{"28}}

\begin{document}
\title{Two-variable Conway polynomial and Cochran's derived invariants}
\author{Sergey A. Melikhov}
\address{Steklov Mathematical Institute of Russian Academy of Sciences, Moscow, Russia}
\email{melikhov@mi-ras.ru}

\begin{abstract}
We note that the Conway potential function $\Omega_L$ of an $m$-component link $L$, $m>1$, can be expressed as 
$\Omega_L(x_1,\dots,x_m)=\Theta_L\big(\varnabla_L(x_1-x_1^{-1},\dots,x_m-x_m^{-1})\big)$ for a unique 
$\varnabla_L\in\Z[z_1,\dots,z_m]$, where $\Theta_L$ is a certain endomorphism of the additive group of $\Z[x_1^{\pm1},\dots,x_m^{\pm1}]$ 
which depends only on the pairwise linking numbers of the components of $L$.
Motivated by applications to topological isotopy, we study the formal power series $\bar\varnabla_L$, obtained by dividing $\varnabla_L$ 
by the Conway polynomials of the components of $L$.

For a $2$-component link with $\lk(L)=0$, the coefficient $\alpha_{1,\,2k-1}$ of $\bar\varnabla_L(u,v)$ at $uv^{2k-1}$ equals Cochran's 
derived invariant $(-1)^{k+1}\beta^k(L)$. 
While this can be deduced from a result of G.-T. Jin, which he proved using the surgical view of the Alexander polynomial, 
we provide an alternative proof, using Seifert matrices.
Our main result is a formula for the same coefficient $\alpha_{1,2k-1}$ in the geometrically subtler case $\lk(L)=1$.
Namely we express it in terms of generalized Cochran invariants $\beta_F^{ij}(P,Q)$, which were studied by Gilmer--Livingston (when $P=Q$) 
and by Tsukamoto--Yasuhara (when $j=0$) and are closely related to the Cochran pairing in the infinite cyclic covering of a knot.

\end{abstract}
\maketitle

\section{Introduction}\label{intro}

\subsection{Conway potential function} \label{intro-conway}
By a {\it link} we will mean a PL link in $S^3$.
For an $m$-component link $L$ let $\Omega_L(x_1,\dots,x_m)$ denote its {\it Conway potential function}.
There is an axiomatic description of $\Omega_L$ in terms of skein relations (see \cite{Jia}) as well as several explicit constructions 
of $\Omega_L$:
\begin{itemize} 
\item in terms of Seifert surfaces that intersect in clasps \cite{Kau}, \cite{Coo2}, \cite{Ci0}, \cite{DMO} (see \S\ref{cpf});
\item in terms of the Fox calculus applied to the Wirtinger presentation \cite{Ha};
\item in terms of the Fox calculus applied to the Neuwirth presentation \cite{BC}*{\S2};
\item in terms of the sign-refined Reidemeister torsion \cite{Tu0}, \cite{BC}*{\S4};
\item as a normalized Euler characteristic of the Ozsv\'ath--Szab\'o link homology \cite{BC}*{\S3};
\item as a quantum invariant (several constructions; see \cite{Vi}, \cite{Harp} and references there);
\item as a diagrammatic state sum \cite{Sa'}.
\end{itemize}
We have $\Omega_L(x_1,\dots,x_m)\in\Z[x_1^{\pm1},\dots,x_m^{\pm1}]$ for $m>1$ and $(x-x^{-1})\Omega_K(x)\in\Z[x^{\pm1}]$
for a knot $K$. 
If the components of $L$ are colored in $n$ colors according to a coloring function $\chi\:\{1,\dots,m\}\to\{1,\dots,n\}$, 
then $\Lambda:=(L,\chi)$ is called a {\it colored link}, and $\Omega_\Lambda(x_1,\dots,x_n)$ denotes $\Omega_L(x_{\chi(1)},\dots,x_{\chi(m)})$.
Thus we have precisely one variable for each color.

In the one-variable case $\Omega_\Lambda(x)$ can be expressed in the form $(x-x^{-1})^{-1}\nabla_\Lambda(x-x^{-1})$, 
where $\nabla_\Lambda(z)\in\Z[z]$ is called the {\it Conway polynomial}.
This follows immediately%
\footnote{If $\Omega\in\Z[x^{\pm1}]$ satisfies $\Omega(x)=\Omega(-x^{-1})$, then its terms come in pairs $a_kx^k+(-1)^ka_kx^{-k}$, 
which can be expressed as $a_k(x-x^{-1})^k$ up to terms of lower (in absolute value) degrees.}
from the symmetry $\Omega_\Lambda(x)=\Omega_\Lambda(-x^{-1})$.
The Conway polynomial is very easy to work with; it is fully determined by two simple axioms 
($\nabla_{\includegraphics[width=0.4cm]{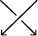}}-\nabla_{\includegraphics[width=0.4cm]{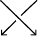}}=z\,\nabla_{\includegraphics[width=0.4cm]{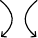}}$
and $\nabla_{\includegraphics[height=0.4cm]{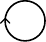}}=1$)
which imply in particular that its coefficients are finite type invariants. 
See e.g.\ \cite{M24-1}*{\S\ref{fti:conway}}.

In the case of $n>1$ variables, $\Omega_\Lambda(x_1,\dots,x_n)$ still satisfies the symmetry 
\[\Omega_\Lambda(x_1,\dots,x_n)=\Omega_\Lambda(-x_1^{-1},\dots,-x_n^{-1})\tag{$\Sigma$}\label{symmetry}\]
(see Lemma \ref{4.1}(a)) but cannot be expressed%
\footnote{For instance, if $U$ denotes the unknot and $U^n$ its pushoff such that $\lk(U,U^n)=n$,
then $\Omega_{(U,U^{-2})}(x,y)=xy+x^{-1}y^{-1}$, which is not a polynomial in $x-x^{-1}$ and $y-y^{-1}$ \cite{Ki}*{proof of Theorem 2}.
Neither is $\Omega_{(U,U^2)}(x,y)=xy^{-1}+x^{-1}y$, since $(xy+x^{-1}y^{-1})-(xy^{-1}+x^{-1}y)=(x-x^{-1})(y-y^{-1})$.}
as a polynomial in $x_1-x_1^{-1},\,\dots,\,x_n-x_n^{-1}$ in general.
However, we may treat $z=x-x^{-1}$ as a quadratic equation in $x$, select one of the two roots $x(z)=\frac{z}2\pm\sqrt{1+\frac{z^2}4}$ and
expand the radical according to the formula $(1+t)^r=1+rt+\frac{r(r-1)}2t^2+\dots$.
Then each $x(z_i)$ becomes a formal power series in $z_i$ (with the same choice of the root for each $i$) which is invertible 
(since its free term is $1$), and so can be substituted for $x_i$ in $\Omega_\Lambda(x_1,\dots,x_n)$.
Thus we get $\Omega_\Lambda(x(z_1),\dots,x(z_n))=\mho_\Lambda(z_1,\dots,z_n)$ for some $\mho_\Lambda\in\Q[[z_1,\dots,z_n]]$, as long as $n>1$.
Both choices of the root lead to the same power series $\mho_\Lambda$,%
\footnote{The Galois group of the quadratic equation $z=x-x^{-1}$ acts on its roots by $x\mapsto-x^{-1}$ 
(since the action defined by this formula takes the equation to itself, or alternatively by using Vieta's formula).
Now the assertion follows from the symmetry (\ref{symmetry}).}
and its coefficients are finite type invariants.%
\footnote{This can be proved similarly to \cite{MuH}*{proof of Lemma 3.2}. 
It is shown in \cite{MuH} and \cite{Dy} that the coefficients of the power series $\Omega_\Lambda(e^{h_1/2},\dots,e^{h_n/2})$ 
are finite type invariants.}

The $2$-variable power series $\mho_\Lambda$ appears in the literature already in the early 90s \cite{HK}, \cite{Kai} under the name 
``two-variable Conway polynomial''.
We refrain from using this terminology, not only because $\mho_\Lambda$ is not a polynomial (to be honest), but also because there is 
another object which better deserves this name, as we will see in a moment.

\subsection{Multi-variable Conway polynomial}
A different approach is found in a paper by M. Kidwell, who observed that the two-variable potential function can be expressed in the form 
\[\Omega_\Lambda(x,y)=\tilde P(u,v)+\bar w\tilde Q(u,v),\tag{$\tilde*$}\label{kidwell}\]
where $u=x-x^{-1}$, $v=y-y^{-1}$ and $\bar w=xy+x^{-1}y^{-1}$, for unique $\tilde P,\tilde Q\in\Z[u,v]$ \cite{Ki}*{Theorem~2}.
This is a purely algebraic observation, which applies to any Laurent polynomial $\Omega$ satisfying the symmetry $\Omega(x,y)=\Omega(-x^{-1},-y^{-1})$.

Our first observation is that the following modified version of Kidwell's decomposition:
\[\Omega_\Lambda(x,y)=P(u,v)+\tfrac{w}2\,Q(u,v),\tag{$*$}\label{kidwell'}\]
where $w=xy^{-1}+x^{-1}y$, extends neatly to any number of variables (Theorem \ref{decomposition}).
The two decompositions are equivalent since $\bar w-w=uv$.
In fact, $w\bar w=u^2+v^2+4$, so by Vieta's theorem $w$ and $-\bar w$ are the roots of the equation $w^2+uvw-(u^2+v^2+4)=0$.
In particular, the two-variable $\mho_\Lambda$ can be obtained from the right hand side of (\ref{kidwell'}) by expanding the radical in 
the formula $w=-\frac{uv}2+\sqrt{\frac{u^2v^2}4+u^2+v^2+4}$.

Our second observation is that no monomial $u^iv^j$ can enter both $P$ and $Q$ with nonzero coefficients, and which of $P$ and $Q$ it can enter 
with a nonzero coefficient is determined by the parities of $m_1$, $m_2$ and $\lk(\Lambda)$, where $m_i$ is the number of components
in the $i$-colored sublink $L_i$ of $\Lambda$ and $\lk(\Lambda)=\lk(L_1,L_2)$.%
\footnote{Namely, $u^iv^j$ can enter $P$ with a nonzero coefficient only if $i\equiv\lk(\Lambda)+m_1\pmod 2$ and $j\equiv\lk(\Lambda)+m_2\pmod 2$, 
and it can enter $Q$ with a nonzero coefficient only if $i\not\equiv\lk(\Lambda)+m_1\pmod 2$ and $j\not\equiv\lk(\Lambda)+m_2\pmod 2$
(see Proposition \ref{parities}).}
Consequently, both $P$ and $Q$ can be reconstructed from
\[\varnabla_\Lambda(u,v):=P(u,v)+Q(u,v),\]
as long as the knowledge of $\lk(\Lambda)$ or at least of its parity is assumed.
(We regard the numbers $m_i$ as known once the colored link $\Lambda$ is given.)
When $\Lambda$ is a two-component link, $\varnabla_\Lambda(0,0)=\Omega_\Lambda(1,1)$ equals $\lk(\Lambda)$, so in this case $\varnabla_\Lambda$ 
contains precisely the same information as $\Omega_\Lambda$.
But for a two-colored link $\Lambda$ of more than two components the integer $\lk(\Lambda)$ is not a function of $\Omega_\Lambda$,%
\footnote{For instance, let $\Lambda=H^{11}\#_1 H^{12}\#_2 H^{22}$ and $\Lambda'=H^{21}\#_1 H^{12}\#_2 H^{21}$, where $H^{ij}$ 
denotes the Hopf link with components of colors $i$ and $j$, and $\#_i$ denotes a connected sum along components of color $i$.
Thus each of $\Lambda$ and $\Lambda'$ has two components of color $1$ and two components of color $2$.
We have $\Omega_\Lambda(x,y)=\Omega_{\Lambda'}(x,y)=(x-x^{-1})(y-y^{-1})$ but $\lk(\Lambda)=1$ and $\lk(\Lambda')=3$.}
and its parity is not a function of $\varnabla_\Lambda$,%
\footnote{For instance, let $H_n^{ij}=(U,U^n)$, where $U$ is the unknot of color $i$ and $U^n$ is its pushoff of color $j$ 
such that $\lk(U,U^n)=n$ and let $\Lambda=H_2^{12}\#_1 H_1^{11}$ and $\Lambda'=H_2^{11}\#_1 H_1^{12}$.
Thus each of $\Lambda$ and $\Lambda'$ has two components of color $1$ and one component of color $2$.
We have $\Omega_\Lambda(x,y)=(xy^{-1}+x^{-1}y)(x-x^{-1})$ and $\Omega_{\Lambda'}(x,y)=2(x-x^{-1})$, whence
$\nabla_\Lambda(u,v)=\nabla_{\Lambda'}(u,v)=2u$, but $\lk(\Lambda)=2$ and $\lk(\Lambda')=1$.}
so in this case we can only say that the pair $\big(\varnabla_\Lambda,\lk(\Lambda)\big)$ contains precisely the same information as 
the pair $\big(\Omega_\Lambda,\lk(\Lambda)\big)$.
In \S\ref{mvc} we extend these observations to any number of colors (Proposition \ref{parities} and Corollary \ref{information}).

We refer to \S\ref{mvc} for the definition of $\varnabla_\Lambda$ in the above spirit for any number of variables,
and here we will be content with stating its less explicit, but more invariant description:

\begin{mainthm}\label{main1}
The Conway potential function $\Omega_\Lambda$ of an $m$-component link $\Lambda$ colored in $n$ colors, where $m>1$, 
can be expressed in the form 
\[\Omega_\Lambda(x_1,\dots,x_n)=\Theta_\Lambda\big(\varnabla_\Lambda(x_1-x_1^{-1},\dots,x_n-x_n^{-1})\big)\] 
for a unique $\varnabla_\Lambda\in\Z[z_1,\dots,z_n]$, where $\Theta_\Lambda$ is a certain endomorphism of the additive group of 
$\Z[x_1^{\pm1},\dots,x_m^{\pm1}]$ which depends only on the pairwise linking numbers of the components of $\Lambda$ and on
the number of components of $\Lambda$ of each color.
\end{mainthm}

\begin{remark} Since $xy^{-1}+x^{-1}y=2+\dfrac{(x-y)^2}{xy}$, we have 
\[\Omega_\Lambda(x,y)=\varnabla_\Lambda(u,v)+\frac{(x-y)^2}{2xy}Q(u,v).\]
Here the remainder $\frac{(x-y)^2}{2xy}Q(u,v)$ is redundant, in the sense that,
modulo the knowledge of $\lk(\Lambda)$, all the information that it contains is already contained in 
$\varnabla_\Lambda(u,v)$.
Yet it is only this redundant remainder that is responsible for the facts that $\mho_\Lambda(u,v)$
is not a polynomial (nor even a rational power series) and has fractional coefficients.
\end{remark}

In \S\ref{reduced} we note that the results of \cite{M24-1} easily imply the following

\begin{mainthm} \label{extension}
Let $L$ be an $m$-component link, with components $K_1,\dots,K_m$.

(a) Each coefficient of the formal power series
\[\bar\varnabla_L(z_1,\dots,z_m):=\dfrac{\varnabla_L(z_1,\dots,z_m)}{\nabla_{K_1}(z_1)\cdots\nabla_{K_m}(z_m)}\]
assumes the same value on all sufficiently close $C^0$-approximations of any given topological link. 
The resulting extension of $\bar\varnabla_L$ by continuity to topological links is an invariant of (non-ambient) isotopy.

(b) Moreover, for each $r$ the extension of $\dfrac{\partial^r\bar\varnabla_L}{\partial z_i^r}(z_1,\dots,z_m)\Big|_{z_i=0}$
to topological links is invariant under sufficiently small $C^0$-perturbation of the $i^{\text{th}}$ component.
\end{mainthm}

\begin{remark} \label{bingsling}
Let us note that exactly how close is ``sufficiently close'' is determined individually for each coefficient.
So the extension of $\bar\varnabla_L$ to topological links need not be a rational power series.
But whenever it is not rational for some topological link $L$, we immediately know that $L$ is not isotopic to any PL link.
This property of $\bar\varnabla_L$ has no counterparts for the rational function 
$\bar\Omega_L:=\Omega_L(x_1,\dots,x_m)/\big(\nabla_{K_1}(x_1-x_1^{-1})\cdots\nabla_{K_m}(x_m-x_m^{-1})\big)$ and for the power series 
$\bar\mho_L:=\mho_L(z_1,\dots,z_m)/\big(\nabla_{K_1}(z_1)\cdots\nabla_{K_m}(z_m))$, and it makes $\bar\varnabla_L$ a powerful invariant 
of isotopy of topological links.

Moreover, if the power series $\dfrac{\partial^r\bar\varnabla_{L}}{\partial u^r}(u,v)\Big|_{u=0}$ is not rational for some 
$2$-component topological link $L$, then we immediately know that $L$ is not isotopic to any link whose second component is PL.

These properties of $\bar\varnabla$ are applied in \cite{M24-3} to obtain some partial results on Rolfsen's 1974 problem: 
{\it Is every topological knot isotopic to a PL knot? In particular, is the Bing sling isotopic to a PL knot?}
\end{remark}

\subsection{Coefficients of the reduced two-variable Conway polynomial}
In the present paper we deal only with PL links (apart from \S\ref{reduced}), but it turns out that for them $\bar\varnabla_L$ 
is of some interest already.
Let $L$ be a $2$-component link, and let us write
\[\bar\varnabla_L(u,v)=\sum_{i=0}^\infty\sum_{j=0}^\infty\alpha_{ij}(L)u^iv^j.\]

\begin{mainthm} \label{main-jin} When $\lk(L)=0$, each $\alpha_{1,2k-1}(L)$ is Cochran's derived invariant invariant $(-1)^{k+1}\beta^k(L)$. 
\end{mainthm} 

The equality of Theorem \ref{main-jin} can be deduced up to a sign from a result of G.-T. Jin \cite{Jin}*{Theorem 4},
if we take into account the equivalence of Cochran's invariants and the Kojima $\eta$-function \cite{Co1}*{Theorem 7.1} 
(see also Theorem \ref{cochran expansion} for a generalization) and Bailey's theorem (see \cite{Jin}*{Lemma 8} or
Theorem \ref{jin}(a)).
The sign can be verified by a direct computation (cf.\ \cite{M04}).

However, Jin's proof, based on the surgical view of the Alexander polynomial, seems to be quite obscure for the purposes of
understanding geometric meaning of link invariants, as it consists in performing a lot of surgeries on $S^3$ and it is difficult 
to keep track of where various objects go under a composition of so many surgeries.
Thus we include an entirely independent proof of Theorem \ref{main-jin}, in terms of Seifert matrices (see Theorem \ref{jin}(b)).

Let us note that all previously known definitions of Cochran's invariants do not extend to the case $\lk(L)=1$.
This is somewhat typical of link theory: invariants whose geometric properties are well-understood often fail to say anything at all
about links of linking number $1$.
This failure seems to be not accidental, but due to a different, and often more complicated, geometry in the case $\lk=1$ (as compared to $\lk=0$).
However, the case $\lk=1$ is crucial in approaching Rolfsen's problem (see Remark \ref{bingsling}).

\begin{lemma} \cite{AADG} \label{C-complex}
Let $L=(K_1,K_2)$ be a link with $\lk(L)=1$.
Then there exist Seifert surfaces $F_i$ for $K_i$, intersecting transversely along a single {\rm clasp arc},
that is, so that $F_1\cap F_2$ is a single arc with one endpoint in $K_1$ and another in $K_2$.
\end{lemma}

A shorter exposition of the same proof from \cite{AADG} can be found in \cite{M24-3}*{Lemma \ref{rolf:C-complex}}.

The following is the main result of the present paper (see Theorem \ref{factorization2}).
It is applied in the paper \cite{M24-3} to obtain its deepest results on Rolfsen's problem.

\begin{mainthm} \label{mainmain} Let $L=(K_1,K_2)$ be a link with $\lk(L)=1$. Then each 
\[\alpha_{1,2n-1}(L)=(-1)^{n+g_1}\sum_{i=1}^{g_1}\beta_{F_2}^{n-1,\,n}(A_{2i},A_{2i-1}),\]
where $F_1$, $F_2$ are Seifert surfaces for $K_1$ and $K_2$, intersecting transversely along a single clasp arc,
and $A_1,\dots,A_{2g_1}$ is a symplectic basis for $F_1$, disjoint from the clasp arc.
\end{mainthm}

Here $\beta_F^{ij}(P,Q)$ are a generalization of Cochran derived invariants, which were studied by Gilmer--Livingston (when $P=Q$) 
and by Tsukamoto--Yasuhara (when $j=0$).
They are defined in \S\ref{generalized cochran}.
Let us note that $\beta_F^{ij}(P,Q)$ are not link invariants in general, as they may depend of the surface $F$.
However, they may be regarded as an alternative (geometric) description of the Cochran pairing \cite{Co1}*{\S7}
in the infinite cyclic covering of a knot, which is in turn a refinement of the Blanchfield pairing.

\S\S\ref{cochran-pairing-section}--\ref{generalized cochran} are devoted to auxiliary results which are used 
in the proofs of Theorems \ref{main-jin} and \ref{mainmain}.
These sections are closely related to the existing literature, but we take this opportunity to fill in some gaps 
and provide missing details.
Thus \S\ref{cochran-pairing-section} contains a proof of Cochran's assertions on the properties of his pairing
(which does not seem to appear in the literature); \S\ref{przytycki-yasuhara} contains an exposition of 
the Przytycki--Yasuhara computation of the Cochran pairing (hopefully more readable than the original version);
and \S\ref{generalized cochran} is mostly concerned with generalization and symmetrization (with respect to
the involution $t\mapsto t^{-1}$ in $\Z[t^{\pm1}]$) of some results by Tsukamoto--Yasuhara and Gilmer--Livingston.

\section{Multi-variable Conway polynomial}\label{mvc}

The only properties of $\Omega_L$ to be used in the present section are as follows.

\begin{lemma} \label{4.1} {\rm (compare \cite{Con}, \cite{Ki})}
Let $L$ be a colored link with $m$ components.
 
(a) $\Omega_L(x_1,\dots,x_n)=\Omega_L(-x_1^{-1},\dots,-x_n^{-1})$.

(b) For $m>1$ the total degree of every nonzero term of $\Omega_L$ is congruent $\bmod 2$ to $m$.

(c) For $m>1$ the exponent of $x_i$ in every nonzero term of $\Omega_L$ is congruent $\bmod 2$ 
to $l_i+m_i$, where $m_i$ is the number of components of color $i$ and $l_i=\lk(L_i,\,\Lambda\but L_i)$, 
where $L_i$ is the sublink of color $i$.
\end{lemma}

The following proof is based on Hartley's definition of $\Omega_L$ \cite{Ha}, which we do not review here.
Later on we give an alternative proof of Lemma \ref{4.1} (see Corollary \ref{4.1'}), based on Cimasoni's construction of $\Omega_L$, 
which we do review in detail in \S\ref{cpf}.

\begin{proof}[Proof. (c)]
It suffices to prove the assertion in the case where all the components of $L$ have distinct colors.
In this case Hartley defines $\Omega_L(x_1,\dots,x_m)$ as a normalized version of his sign-refined Alexander polynomial $\Delta_L(t_1,\dots,t_m)$, 
which is well-defined up to multiplication by monomials $t_1^{i_1}\dots t_m^{i_m}$ \cite{Ha}.
Namely, $(x-x^{-1})\Omega_K(x)=x^{\mu}\Delta_K(x^2)$ for a knot $K$, and for a link $L$ with $m>1$ components, 
$\Omega_L(x_1,\dots,x_m)=x_1^{\mu_1}\dots x_n^{\mu_m}\Delta_L(x_1^2,\dots,x_m^2)$, where the integers $\mu,\mu_1,\dots,\mu_n$ 
are uniquely determined by the symmetry relation $\Omega_L(x_1,\dots,x_m)=(-1)^m\Omega_L(x_1^{-1},\dots,x_m^{-1})$
of \cite{Ha}*{(5.5)}.

The desired assertion in the case where all the components of $L$ have distinct colors follows immediately from the fact, 
noted without proof in \cite{Tr2}*{\S7}, that the parity of $\mu_i$ is opposite to that of $l_i$.
For the proof of this fact we observe that the definition of $\mu_i$ in \cite{Ha}*{(2.4)}, which depends on the choice of 
a plane diagram of $L$, implies that $\mu_i$ has the same parity as $l_i+\delta_i+\sigma_i$, where $\sigma_i$ is the number
of positively oriented Seifert circles minus the number of negatively oriented Seifert circles of the plane diagram $D_i$
of the $i$th component of $L$ and $\delta_i$ is the number of double points of the plane curve $D_i$.
Now the parity of $\delta_i$ is opposite to that of the turning number of $D_i$ \cite{Wh}*{Theorem 2}, and the turning number
of $D_i$ is easily seen to be equal to $\sigma_i$.
\end{proof}

\begin{proof}[(b)] This follows easily from (c).
\end{proof}

\begin{proof}[(a)]
For links with $>1$ components this follows from (b) and from the symmetry relation 
$\Omega_L(x_1,\dots,x_n)=(-1)^m\Omega_L(x_1^{-1},\dots,x_n^{-1})$, which is a consequence of \cite{Ha}*{(5.5)}.
For knots the assertion follows from $\Omega_K(x)=\nabla_K(x-x^{-1})/(x-x^{-1})$.
\end{proof}

\begin{theorem} \label{decomposition}
The Conway potential function of any colored link $L$ of $m>1$ components can be uniquely written in the form
\[\Omega_L(x_1,\dots,x_n)=\frac12\hspace{-17pt}\sum_{\substack{k\ge 0\\ 1\le i_1<\dots<i_{2k}\le n\vphantom{A^E}}}\hspace{-20pt}
\left(\frac{x_{i_1}}{x_{i_2}}\cdots\frac{x_{i_{2k-1}}}{x_{i_{2k}}}+\frac{x_{i_2}}{x_{i_1}}\cdots\frac{x_{i_{2k}}}{x_{i_{2k-1}}}\right)
P_{i_1,\dots,i_{2k}}(x_1-x_1^{-1},\,\dots,\,x_n-x_n^{-1})\]
for some polynomials $P_{i_1,\dots,i_{2k}}$ with integer coefficients.
Moreover, when $n=2l$, the coefficients of the polynomial $P_{1,2,\dots,n}$ are even.
\end{theorem}

Following Conway \cite{Con} and Kidwell \cite{Ki}, we will use the abbreviation 
\[\big\{f(x_1,\dots,x_n)\big\}:=f(x_1,\dots,x_n)+f(-x_1^{-1},\dots,-x_n^{-1})\] 
for any function $f(x_1,\dots,x_n)$.
In this notation the formula of Theorem \ref{decomposition} becomes
\[\Omega_L(x_1,\dots,x_n)=\frac12\sum_{\substack{k\ge 0\\ 1\le i_1<\dots<i_{2k}\le n\mathstrut}}
\left\{\frac{x_{i_1}}{x_{i_2}}\cdots\frac{x_{i_{2k-1}}}{x_{i_{2k}}}\right\}
P_{i_1,\dots,i_{2k}}\big(\{x_1\},\,\dots,\,\{x_n\}\big).\tag{$**$}\label{decomposition-formula}\]
Let us note that $\{1\}=2$, which occurs in the summand corresponding to $k=0$ in (\ref{decomposition-formula}).

\begin{remark}
The proof of Theorem \ref{decomposition} works for any Laurent polynomial $\Omega$ satisfying the symmetry
(\ref{symmetry}).
\end{remark}

\begin{remark}
The assertion on integrality of $P_{i_1\dots i_{2k}}$ in
Theorem \ref{decomposition} will not hold already for $n=3$ (respectively $n=5$) if
$\big\{\frac{x_{i_1}}{x_{i_2}}\cdots\frac{x_{i_{2k-1}}}{x_{i_{2k}}}\big\}$ is
replaced with $\{x_{i_1}\cdots x_{i_{2k}}\}$ (respectively with
$\big\{\frac{x_{i_1}\cdots x_{i_k}}{x_{i_{k+1}}\cdots x_{i_{2k}}}\big\}$)
in (\ref{decomposition-formula}), at least for some Laurent polynomial $\Omega$ satisfying (\ref{symmetry}).
\end{remark}

We will also use the abbreviation \[\big[f(x_1,\dots,x_n)\big]:=f(x_1,\dots,x_n)-f(-x_1^{-1},\dots,-x_n^{-1})\]
though this one will be needed less frequently. 

\begin{lemma}\label{identities} Let $M$ be a monomial in $x_1^{\pm 1},\dots,x_n^{\pm 1}$.

(a) $\{x_iM\}-\{x_i^{-1}M\}=\{x_i\}\{M\}$.

(b) $[x_iM]-[x_i^{-1}M]=\{x_i\}[M]$.

(c) $\{x_ix_jM\}+\{x_i^{-1}x_j^{-1}M\}=\{x_ix_j\}\{M\}$.

(d) $2\{x_ix_j^{-1}M\}=\{x_i\}\{x_j^{-1}M\}+\{x_j^{-1}\}\{x_iM\}+\{x_ix_j\}\{M\}$. 
\end{lemma}

\begin{proof}[Proof. (a), (b)]
Let $d$ be the total degree of $M$.
When $d$ is even, the two desired identities follow respectively from the following two identities:
\begin{gather*}
(x_iM-x_i^{-1}M^{-1})-(x_i^{-1}M-x_iM^{-1})=(x_i-x_i^{-1})(M+M^{-1})\\
(x_iM+x_i^{-1}M^{-1})-(x_i^{-1}M+x_iM^{-1})=(x_i-x_i^{-1})(M-M^{-1})
\end{gather*}
When $d$ is odd, the same holds if the two desired identities are interchanged.
\end{proof}

\begin{proof}[(c)] Similarly to (a).
\end{proof}

\begin{proof}[(d)] This follows from (a) and (c).
\end{proof}

\begin{proof}[Proof of Theorem \ref{decomposition}] {\it Existence.}
By the symmetry relation (\ref{symmetry}), $\Omega_L$ includes together with every term
$Ax_1^{p_1}\cdots x_n^{p_n}$ the term $(-1)^{p_1+\dots+p_n}Ax_1^{-p_1}\cdots x_n^{-p_n}$,
and so can be written as a $\Z$-linear combination of the Laurent polynomials
$\{x_1^{p_1}\cdots x_n^{p_n}\}$.
Using the formula of Lemma \ref{identities}(a) one can express each $\{x_1^{p_1}\dots x_n^{p_n}\}$ in the form
\[\sum_{\substack{k\ge 0\\ 1\le i_1<\dots<i_k\le n\mathstrut}}\big\{x_{i_1}x_{i_2}^{-1}x_{i_3}x_{i_4}^{-1}\cdots
(x_{i_{k-1}}x_{i_k}^{-1})^{(-1)^k}\big\} P'_{i_1,\dots,i_k}(\{x_1\},\dots,\{x_n\})\]
for some $P'_{i_1,\dots,i_k}\in\Z[z_1,\dots,z_n]$, $k>0$, and some $P'\in\frac12\Z[z_1,\dots,z_n]$.
The summands corresponding to $k=1$ can be included in $P'$, and one can get
rid of the summands corresponding to odd $k\ge 3$ by repeated use of
the formulas (a) and (d) of Lemma \ref{identities}.
But it is not clear from this approach that the resulting polynomials will have
half-integer coefficients.
To see this, represent
$2\big\{x_{i_1}x_{i_2}^{-1}x_{i_3}\cdots x_{i_{2k}}^{-1}x_{i_{2k+1}}\big\}$ as
\begin{multline*}\big\{x_{i_1}x_{i_2}^{-1}x_{i_3}\cdots x_{i_{2k}}^{-1}x_{i_{2k+1}}\big\}-
\big\{x_{i_1}^{-1}x_{i_2}x_{i_3}^{-1}\cdots x_{i_{2k}}x_{i_{2k+1}}^{-1}\big\}\\
+(1-1)\sum_{j=1}^{2k}\big\{x_{i_1}^{-1}x_{i_2}x_{i_3}^{-1}\cdots
(x_{i_{j-2}}x_{i_{j-1}}^{-1}x_{i_j}x_{i_{j+1}}x_{i_{j+2}}^{-1}x_{i_{j+3}})^{(-1)^j}\cdots x_{i_{2k}}^{-1}x_{i_{2k+1}}\big\}.
\end{multline*}
Then the formula of Lemma \ref{identities}(a) yields
\begin{multline*}
2\big\{x_{i_1}x_{i_2}^{-1}x_{i_3}\cdots x_{i_{2k}}^{-1}x_{i_{2k+1}}\big\}\\
=\sum_{j=1}^{2k+1}(-1)^{j+1}\{x_{i_j}\}
\big\{x_{i_1}^{-1}x_{i_2}x_{i_3}^{-1}\cdots
(x_{i_{j-2}}x_{i_{j-1}}^{-1}x_{i_{j+1}}x_{i_{j+2}}^{-1}x_{i_{j+3}})^{(-1)^j}\cdots x_{i_{2k}}^{-1}x_{i_{2k+1}}\big\}.
\end{multline*}
Thus each $P'_{i_1,\dots,i_{2k+1}}$ can be dispensed with at the cost of bringing in half-integer coefficients to the polynomials
$P'_{j_1,\dots,j_{2k}}$.
Note that when $n=2l$, the polynomial $P'_{1,\dots,n}$ is not affected by this process, and so its coefficients remain integer.

{\it Uniqueness.}
It remains to verify the uniqueness of the decomposition (\ref{decomposition-formula}).
Suppose, by way of contradiction, that a nontrivial expression
$Q(x_1,\dots,x_n)$ in the form of the right hand side of (\ref{decomposition-formula}) is
identically zero.
Then so is $Q(x_1,\dots,x_{n-1},x_n)-Q(x_1,\dots,x_{n-1},-x_n^{-1})$, which
can be rewritten as $$[x_n]\sum_{\substack{k\ge 1\\ 1\le i_1<\dots<i_{2k-1}\le n-1\mathstrut}}
\big[x_{i_1}x_{i_2}^{-1}x_{i_3}\cdots x_{i_{2k-2}}^{-1}x_{i_{2k-1}}\big]
P_{i_1,\dots,i_{2k-1},n}(\{x_1\},\dots,\{x_n\})=0.$$
Let us denote the left hand side by $[x_n]R(x_1,\dots,x_n)$; then $R(x_1,\dots,x_n)$
is identically zero.
Hence so is
$R(x_1,\dots,x_{n-2},x_{n-1},x_n)-R(x_1,\dots,x_{n-2},-x_{n-1}^{-1},x_n)$,
which can be rewritten as
$$[x_{n-1}]\sum_{\substack{k\ge 0\\ 1\le i_1<\dots<i_{2k}\le n-2\mathstrut}}
\big\{x_{i_1}x_{i_2}^{-1}\cdots x_{i_{2k-1}}x_{i_{2k}}^{-1}\big\}
P_{i_1,\dots,i_{2k},n-1,n}(\{x_1\},\dots,\{x_n\})=0.$$
Repeating this two-step procedure $\lfloor\frac n2\rfloor$ times, we will
end up with $$[x_1]\{1\}P_{1,\dots,n}(\{x_1\},\dots,\{x_n\})=0
\quad\text{ or }\quad
[x_2]\{1\}P_{2,\dots,n}(\{x_1\},\dots,\{x_n\})=0$$
according as $n$ is even or odd.
Consider, for example, the case of odd $n$.
In this case we conclude that $P_{2,\dots,n}=0$.
But then by symmetry $P_{1,\dots,\hat\imath,\dots,n}=0$ for each $i$.
Returning to the previous stage
\[[x_4]\sum_{\substack{k\ge 0\\ 1\le i_1<\dots<i_{2k}\le 3\mathstrut}}
\big\{x_{i_1}x_{i_2}^{-1}\cdots x_{i_{2k-1}}x_{i_{2k}}^{-1}\big\}
P_{i_1,\dots,i_{2k},4\dots n}(\{x_1\},\dots,\{x_n\})=0,\] we can now
substitute zeroes for $P_{2,3,4,\dots,n},P_{1,3,4,\dots,n},P_{1,2,4,\dots,n}$,
and so we get $P_{4,\dots,n}=0$.
Continuing to the earlier stages, we will similarly verify that each
$P_{i_1,\dots,i_{2k}}=0$.
\end{proof}

Lemma \ref{4.1}(b,c) implies

\begin{proposition}\label{parities} 
In the notation of Theorem \ref{decomposition}, every nonzero term $T$ of each 
polynomial $P_{i_1,\dots,i_{2k}}(z_1,\dots,z_n)$ satisfies:

(a) the total degree of $T$ is congruent $\bmod 2$ to $m$;

(b) the exponent of $z_i$ in $T$ is congruent $\bmod 2$ to $m_i+l_i+\chi_i$, where
\begin{itemize}
\item $m_i$ is the number of components of color $i$;
\item $l_i$ is the linking number between the sublink of color $i$ and the remaining sublink;
\item $\chi_i=1$ if $i\in\{i_1,\dots,i_{2k}\}$, and $\chi_i=0$ otherwise.
\end{itemize}
\end{proposition}

\begin{example} \label{parities-example} Let $L$ be a two-component link.
If $\lk(L)$ is odd, then $P(u,v)$ contains only monomials of the form $u^{2i}v^{2j}$
and $P_{1,2}(u,v)$ contains only monomials of the form $u^{2i+1}v^{2j+1}$.
If $\lk(L)$ is even, then $P(u,v)$ contains only monomials of the form $u^{2i+1}v^{2j+1}$ 
and $P_{1,2}(u,v)$ contains only monomials of the form $u^{2i}v^{2j}$.
\end{example}

In general, Proposition \ref{parities} implies that if $c_{e_1,\dots,e_n} z_1^{e_1}\cdots z_n^{e_n}$ is 
a nonzero term of $P_{i_1,\dots,i_{2k}}$, then the indices $i_1,\dots,i_{2n}$ are fully determined by 
the parities of the exponents $e_1,\dots,e_n$.
Consequently the individual polynomials $P_{i_1,\dots,i_{2k}}$ can be reconstructed from their sum%
\footnote{This definition of $\varnabla_L$ makes sense when $L$ has at least two components. 
For a knot $K$ we set $\varnabla_K(z)=z^{-1}\nabla_K(z)$, where $\nabla_K(z)$ is the Conway polynomial.}
\[\varnabla_L(z_1,\dots,z_n):=\sum_{\substack{k\ge 0\\ 1\le i_1<\dots<i_{2k}\le n\mathstrut}}P_{i_1,\dots,i_{2k}}(z_1,\,\dots,\,z_n).\]
This reconstruction of course depends on Proposition \ref{parities}, which in turn involves knowledge of the numbers $m_1,\dots,m_n$ 
and $l_1,\dots,l_n$.
We regard the $m_i$ as given when we are given the colored link $L$; but the linking numbers $l_i$ strictly speaking
require computation.
Thus we get

\begin{corollary}\label{information} 
Assuming that the linking numbers $l_i$ are given, $\varnabla_L$ contains 
the same information as $\Omega_L$.
\end{corollary}

Theorem \ref{main1} is an abridged form of Theorem \ref{main1'} and Remark \ref{endomorphism}.

\begin{theorem}\label{main1'}
The Conway potential function $\Omega_\Lambda$ of an $m$-component link $\Lambda$ colored in $n$ colors, where $m>1$, 
can be expressed in the form 
\[\Omega_\Lambda(x_1,\dots,x_n)=\Theta_\Lambda\big(\varnabla_\Lambda(z_1,\dots,z_n)\big)\] 
for a unique $\varnabla_\Lambda\in\Z[z_1,\dots,z_n]$, where $\Theta_\Lambda\:\Z[z_1,\dots,z_n]\to\Z[x_1^{\pm1},\dots,x_n^{\pm1}]$ is 
a homomorphism of the additive groups which depends only on the number of components $m_i$ in the sublink $L_i$ of 
the $i^{\text{th}}$ color, for each $i$, and on the linking numbers $l_i:=\lk(L_i,\,\Lambda\but L_i)$.

Namely, $\Theta_\Lambda$ is given on the additive generators by $\Theta_\Lambda(z_1^{e_1}\cdots z_n^{e_n})=0$ if 
$e_1+\dots+e_n\not\equiv m\pmod 2$, and otherwise by
\[\Theta_\Lambda(z_1^{e_1}\cdots z_n^{e_n})=
\frac12\left(\frac{x_{i_1}}{x_{i_2}}\cdots\frac{x_{i_{2k-1}}}{x_{i_{2k}}}+\frac{x_{i_2}}{x_{i_1}}\cdots\frac{x_{i_{2k}}}{x_{i_{2k-1}}}\right)
(x_1-x_1^{-1})^{e_1}\cdots(x_n-x_n^{-1})^{e_n},\]
where $i_1<\dots<i_{2k}$ are uniquely determined by $i\in\{i_1,\dots,i_{2k}\}\Leftrightarrow e_i\not\equiv l_i+m_i\pmod2$. 
\end{theorem}

\begin{remark} The congruences $e_i\equiv l_i+m_i\pmod2$ fail for an even number of indices $i$, since their sum
$e_1+\dots+e_n\equiv l_1+\dots+l_n+m_1+\dots+m_n$ always holds (due to $l_1+\dots+l_m=2\sum_{1\le i<j\le m}\lk(K_i,K_j)$, where
$K_1,\dots,K_m$ are the components of $\Lambda$).
\end{remark}

\begin{remark} \label{endomorphism}
If each $z_i$ is identified with $x_i-x_i^{-1}$, then $\Theta_\Lambda$ is easily seen to extend to an endomorphism 
$\bar\Theta_\Lambda$ of the additive group of $\Z[x_1^{\pm1},\dots,x_m^{\pm1}]$. 
Namely, $\bar\Theta_\Lambda$ is given on the additive generators by $\bar\Theta_\Lambda(x_1^{e_1}\cdots x_n^{e_n})=0$ if 
$e_1+\dots+e_n\not\equiv m\pmod 2$, and otherwise by
\[\bar\Theta_\Lambda(x_1^{e_1}\cdots x_n^{e_n})=
\frac12\left(\frac{x_{i_1}}{x_{i_2}}\cdots\frac{x_{i_{2k-1}}}{x_{i_{2k}}}+\frac{x_{i_2}}{x_{i_1}}\cdots\frac{x_{i_{2k}}}{x_{i_{2k-1}}}\right)
x_1^{e_1}\cdots x_n^{e_n},\]
where $i_1<\dots<i_{2k}$ are uniquely determined by $i\in\{i_1,\dots,i_{2k}\}\Leftrightarrow e_i\not\equiv l_i+m_i\pmod2$. 
\end{remark}

\begin{proof}[Proof of Theorem \ref{main1'}]
Let $\varnabla_L$ be defined as above (based on the existence part of Theorem \ref{decomposition}).
Then Proposition \ref{parities} implies that $\Omega_L=\Theta_L(\varnabla_L)$; and the uniqueness part of 
Theorem \ref{decomposition} implies that the latter relation uniquely determines $\varnabla_L$.
\end{proof}

\section{Reduced multi-variable Conway polynomial} \label{reduced}
One apparent downside of $\varnabla_\Lambda$ as compared with $\mho_\Lambda$ is that the coefficients of $\varnabla_\Lambda$ are not
finite type invariants in general. 
However, since $\Omega_\Lambda=\Omega_\Lambda(x_1,\dots,x_n)$ satisfies the skein relation
\[\Omega_{\includegraphics[width=0.5cm]{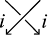}}-\Omega_{\includegraphics[width=0.5cm]{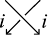}}=(x_i-x_i^{-1})\,\Omega_{\includegraphics[width=0.5cm]{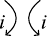}},\]
$\varnabla_\Lambda=\varnabla_\Lambda(z_1,\dots,z_n)$ satisfies the skein relation
\[\varnabla_{\includegraphics[width=0.5cm]{1+.pdf}}-\varnabla_{\includegraphics[width=0.5cm]{1-.pdf}}=z_i\,\varnabla_{\includegraphics[width=0.5cm]{1o.pdf}},\]
which together with the vanishing of the coefficients of $\varnabla_\Lambda$ in degrees $<-1$ implies in the usual way
that the coefficients of $\varnabla_\Lambda$ are {\it colored} finite type invariants.

{\it Colored finite type invariants} of links were introduced by Kirk and Livingston \cite{KL} (they did not use the word ``colored'').
An invariant $v$ of links colored in $n$ colors is said to be of {\it type $(k_1,\dots,k_n)$} if its standard extension to singular links 
(see \cite{M24-1}*{\S\ref{fti:fti}}) vanishes on all singular links with $k_1+1$ self-intersections and/or intersections between components of color $1$;
on all singular links with $k_2+1$ self-intersections and/or intersections between components of color $2$; and so on.
(Note that there is no condition on double points that are intersections between components of distinct colors.)
An invariant is of {\it colored finite type} if and only it is of type $(k_1,\dots,k_n)$ for some $k_1,\dots,k_n$ 
(this can be taken as a definition, and to see that it is equivalent to the definition of \cite{KL} is an easy exercise, cf.\ \cite{M24-1}*{\S4}).
The above skein relation for $\varnabla_\Lambda$ implies in the usual way (see \cite{BN}*{proof of Theorem 2}
or \cite{M24-1}*{proof of Lemma \ref{fti:conway-coefficients}}) that the coefficient of $\varnabla_\Lambda(z_1,\dots,z_n)$ at 
$z_1^{k_1}\cdots z_n^{k_n}$ is of type $(k_1,\dots,k_n)$.

What is good about colored finite type invariants is that they satisfy the following property.

\begin{theorem} \cite{M24-1}*{Corollary \ref{fti:main1''}} \label{top-extension}
Let $v$ be a colored finite type invariant of links, where all components are colored in distinct colors.
Suppose that $v$ is invariant under (non-ambient) PL isotopy.
Then $v$ assumes the same value on all sufficiently close $C^0$-approximations of any given topological link. 
Moreover, the extension of $v$ by continuity to topological links is an invariant of isotopy.
\end{theorem}

Where do we get invariants of PL isotopy? Since two links are PL isotopic if and only if they are equivalent under 
the equivalence relation generated by ambient isotopy and insertion of local knots (see \cite{Ro2}*{Theorem 4.2}),
the connected sum formula for $\Omega_L$ (see \cite{Ci}*{Corollary 3.5}) implies that the rational function
\[\bar\Omega_L:=\frac{\Omega_L(x_1,\dots,x_m)}{\nabla_{K_1}(x_1-x_1^{-1})\cdots\nabla_{K_m}(x_m-x_m^{-1})},\]
where $K_1,\dots,K_m$ are the components of $L$, is invariant under PL isotopy (compare \cite{Tr2}*{Theorem 1 and \S5}, \cite{Ro4}).
Moreover, since $\nabla_K(z)$ contains only terms of even degrees for a knot $K$ (see e.g.\ 
\cite{M24-1}*{Lemma \ref{fti:lickorish}(b)}), multiplication by
$\nabla_K(x_i-x_i^{-1})$ commutes with the endomorphism $\bar\Theta_L$ (see Remark \ref{endomorphism}), and it follows that 
the rational function
\[\bar\varnabla_L(z_1,\dots,z_m):=\dfrac{\varnabla_L(z_1,\dots,z_m)}{\nabla_{K_1}(z_1)\cdots\nabla_{K_m}(z_m)}\]
is also invariant under PL isotopy. 
Now $\varnabla_L$ is a genuine polynomial (for $m>1$) and each $\nabla_{K_i}(0)=1$, so we may regard 
$\bar\varnabla_L(z_1,\dots,z_m)$ as a formal power series: $\bar\varnabla_L\in\Z[[z_1,\dots,z_m]]$.
Each coefficient of this power series is invariant under PL isotopy and is a colored finite type invariant
(as a polynomial in colored finite type invariants, see \cite{M24-1}*{Corollary \ref{fti:product}}).
Hence from Theorem \ref{top-extension} we get

\begin{corollary} Each coefficient of $\bar\varnabla_L(z_1,\dots,z_m)$ assumes the same value on all sufficiently 
close $C^0$-approximations of any given topological link. 
The resulting extension of $\bar\varnabla_L$ by continuity to topological links is an invariant of isotopy.
\end{corollary}

In fact, one can say a bit more.
Given a colored link $\Lambda=(L,\chi)$, we can define
$\bar\varnabla_\Lambda(x_1,\dots,x_n)$ in the usual way, as $\bar\varnabla_L(x_{\chi(1)},\dots,x_{\chi(m)})$.
Since the coefficient of $\varnabla_\Lambda(z_1,\dots,z_n)$ at $z_1^{k_1}\cdots z_n^{k_n}$ is of type $(k_1,\dots,k_n)$,
it is easy to see that so is the coefficient of $\bar\varnabla_\Lambda(z_1,\dots,z_n)$ at $z_1^{k_1}\cdots z_n^{k_n}$
(see \cite{M24-1}*{Corollary \ref{fti:product}}).
Hence all the coefficients of $\dfrac{\partial^r\bar\varnabla_\Lambda}{\partial z_i^r}(z_1,\dots,z_n)\Big|_{z_i=0}$
are of types $(k_1,\dots,k_n)$, where $k_i\le r$.
Then by the proof of \cite{M24-1}*{Theorem \ref{fti:main1}} they are all invariant under $r$-quasi-isotopy of $\Lambda$
with support in components of color $i$.
Therefore by the proof of \cite{M24-1}*{Theorem \ref{fti:isotopy0}} their extensions to topological links are all invariant 
under sufficiently close $C^0$-perturbation of components of color $i$, in the following sense:

\begin{theorem} Given a colored topological link $\L$, for each $\epsilon>0$ there exists an $r$ such that
the entire power series \[\dfrac{\partial^r\bar\varnabla_\Lambda}{\partial z_i^r}(z_1,\dots,z_n)\Bigg|_{z_i=0}\]
(where $\bar\varnabla$ is extended by continuity to topological links) assumes the same value on all colored
topological links $\Lambda$ whose components of color $i$ are $\epsilon$-close to those of $\L$, and 
components of the other colors coincide with those of $\L$.
\end{theorem}

We conclude this section with the following remark.

\begin{proposition} \label{parities2}
The total degree of every term of (a) $\varnabla_L$, (b) $\bar\varnabla_L$ has the same parity as number of components of $L$.
\end{proposition}

\begin{proof} Part (a) follows from Proposition \ref{parities}(a).
Part (b) follows from (a) and the fact that the Conway polynomial of a knot contains only terms of even degrees
(see e.g.\ \cite{M24-1}*{Lemma \ref{fti:lickorish}(b)}).
\end{proof}

\section{Cochran pairing} \label{cochran-pairing-section}

Let $K$ be a PL knot in $S^3$.
Let $X=S^3\but K$ and let $p\:\tilde X\to X$ be the infinite cyclic covering.
The group $H_1(X)\simeq\Z\simeq\left<t\mid\,\right>$ acts on $\tilde X$ by covering transformations.
Consequently $H_1(\tilde X)$ is a module over $\Lambda:=\Z[H_1(X)]\simeq\Z[t^{\pm 1}]$.

\begin{lemma} \label{torsion} 
(a) The action by $t-1$ on $H_1(\tilde X)$ and on $H_2(\tilde X)$ is by automorphisms.

(b) $H_2(\tilde X)=0$ and $H_1(\tilde X)$ is $\Lambda$-torsion.
\end{lemma}

We include the proof for convenience.
Part (b) is well-known (see e.g.\ \cite{Kaw}*{Proposition 7.3.4(1)}; see also \cite{Mil}*{Proof of Assertion 5}).

\begin{proof}[Proof. (a)] The short exact sequence $0\to\Lambda\xr{t-1}\Lambda\xr{t=1}\Z\to 0$
of $\Lambda$-modules yields a long exact sequence of cohomology groups of $X$ with local coefficients,%
\footnote{See for instance \cite{Hat}*{\S3.H} concerning homology with local coefficients.}
which upon applying the isomorphism $H_i(X;\,\Lambda)\simeq H_i(\tilde X;\,\Z)$ becomes
\[H_3(X)\to H_2(\tilde X)\xr{t-1}H_2(\tilde X)\to H_2(X)\to H_1(\tilde X)\xr{t-1}H_1(\tilde X)\to H_1(X)\xr{\delta} H_0(\tilde X).\]
We have $H_2(X)=0=H_3(X)$ by the Alexander duality.
The connecting homomorphism $\delta$ sends the class of a $1$-cycle $z$ to the class of a $0$-cycle $w$ 
such that $(t-1)w=\partial\tilde z$, where $p(\tilde z)=z$.
It follows that $\delta$ is an isomorphism.
Therefore the action by $t-1$ is an automorphism both on $H_1(\tilde X)$ and on $H_2(\tilde X)$.
\end{proof}

\begin{proof}[(b)]
Since $H_1(\tilde X)$ is finitely generated over $\Lambda$, which is Noetherian,%
\footnote{$\Z[t]$ is Noetherian (see \cite{DF}*{Theorem 15.3}) and every localization 
of a Noetherian ring is Noetherian (see \cite{DF}*{Proposition 15.38(4)}).
Alternatively, see \cite{Pu2} for a direct proof.}
$H_1(\tilde X)$ is Noetherian.%
\footnote{See \cite{Lang}*{Proposition X.1.4}.}
Hence for each $m\in H_1(\tilde X)$ the submodules $\left<m\right>\subset\big\langle(t-1)^{-1}m\big\rangle
\subset\big\langle(t-1)^{-2}m\big\rangle\subset\dots$ stabilize.
Thus there exist a $k$ and a $p\in\Lambda$ such that $(t-1)^{-k-1}m=p(t-1)^{-k}m$. 
Then $m=p(t-1)m$, whence $\big(p(t-1)-1\big)m=0$.

Thus we have shown that $H_1(\tilde X)$ is $\Lambda$-torsion.
Similarly $H_2(\tilde X)$ is $\Lambda$-torsion.
On the other hand, $X=S^3\but K$ deformation retracts onto a $2$-polyhedron $X'$, and
$\tilde X$ is homotopy equivalent to the infinite cyclic cover $\tilde X'$ of $X'$.
Now $H_2(\tilde X')=\ker\left[C_2(\tilde X')\to C_1(\tilde X')\right]$ is free over $\Z[\Z]$,
since a submodule of a free module over a PID is free.%
\footnote{See \cite{Lang}*{p.\ 880}}
\end{proof}

\begin{remark} A. Zastrow asked me if there is a geometric proof that every $m\in H_1(\tilde X)$
is divisible by $t-1$.
Such a proof is illustrated in \cite{Co1}*{Figure 7.2}. 
Let us discuss some details of this construction, especially since $K$ is assumed to be smooth in \cite{Co1}.
Let $\Sigma$ be an $h$-Seifert surface of $K$ (see \cite{M21}*{Lemma 2.1}) and
let us fix its lift $\tilde\Sigma$ in $\tilde X$.
The general case easily reduces to the case where $m$ is represented by a knot $\tilde Q$ in $\tilde X$ 
which lies in the open region $R$ bounded by $\tilde\Sigma$ and the consecutive lift $t\tilde\Sigma$.
Then $\tilde Q$ is a lift of a knot $Q\subset S^3\but(K\cup\Sigma)$.
Since $\lk(Q,K)=1$, $Q$ bounds a Seifert surface $S$ in $S^3\but K$ (see \cite{M24-3}*{proof of Theorem \ref{rolf:boundary}}).
We may assume that $S$ meets $\Sigma$ transversely.
Let $C=S\cap\Sigma$, and let $\tilde S$ be the closure of the lift of $S\but C$ in $R$.
Then $\partial\tilde S=(t-1)\tilde C$, where $\tilde C$ is a lift of $C$ lying in $\tilde\Sigma$.
\end{remark}

Let us fix some general sign conventions.
For an oriented smooth manifold $N^n$, a tangent frame $(v_1,\dots,v_{n-1})$ of $\partial N$ is {\it positive} 
if $(v_1,\dots,v_n)$ is a positive tangent frame of $N$, where $v_n$ is an inward looking normal vector.
For an oriented codimension one smooth submanifold $F$ of $N$, a normal vector $v_n$ of $F$ is {\it positive}
if $(v_1,\dots,v_n)$ is a positive tangent frame of $N$, where $(v_1,\dots,v_{n-1})$
is a positive tangent frame of $F$.
For oriented smooth submanifolds $L^l$ and $M^m$ of $N$ which meet transversely in a smooth manifold $K^k$ 
(so, $l-k=n-m$), a {\it positive} tangent frame of $K$ is given by the overlap $(v_{l-k+1},\dots,v_l)$ of 
a positive tangent frame $(v_1,\dots,v_l)$ of $L$ and a positive tangent frame $(v_{n-m+1},\dots,v_n)$ 
of $M$ whose combination $(v_1,\dots,v_n)$ is a positive tangent frame of $N$.
When $k=0$, this rule orients individual points of $K$, and then the {\it intersection pairing}
$L\cdot M$ is defined to be their algebraic sum.
If additionally $\partial M=\emptyset$ and $N$ is a sphere, then the {\it linking number}
$\lk(\partial L,M)=L\cdot M$. 

Let us fix some orientations of $S^3$ and of the knot $K$.
When $S^3$ is identified with the one-point compactification of $\R^3$,
we want the standard coordinate frame $(i,j,k)$ of $\R^3$ (where $i$ points right and
$j$ points up in a sheet of paper, and $k$ points toward the reader) to be positive.
Given a Seifert surface $F$ for $K$, its orientation is determined by that of $K$.
Then the orientation of $S^3$ determines a co-orientation of $F$.
The co-orientation of $F$ in turn determines a generator $t$ of the group of covering transformations of $\tilde X$.
This $t$ obviously does not depend on the choice of $F$.
(Alternatively, $t$ could be determined using the Alexander duality $H_1(K)\simeq H^1(X)$ and the isomorphism $H^1(X)\simeq [X,S^1]$, 
which does not involve the choice of a Seifert surface, but involves more steps.)

Let $A$ and $B$ be $1$-cycles in $\tilde X$ such that $A$ and $t^nB$ have disjoint supports for each $n\in\Z$.
Let $f(t)\in\Lambda$ be any annihilator of $[A]\in H_1(\tilde X)$.
Then the $1$-cycle $f(t)A$ bounds a $2$-chain $\zeta$ in $\tilde X$.%
\footnote{An explicit construction of such a $\zeta$, with $f(t)$ being the Alexander polynomial of $K$, can be extracted from 
the proof of Proposition \ref{PY1}(a) below.}
The {\it Cochran pairing} \cite{Co1}*{\S7} is defined by
\[\left<A,B\right>=\frac1{f(t)}\sum_{n=-\infty}^\infty
(\zeta\cdot t^nB)t^n,\]
where $\cdot$ stands for the intersection pairing in the oriented $3$-manifold $\tilde X$, and only finitely many 
of the summands may be nonzero.
Thus $\left<A,B\right>\in\Q(t)$, the field of fractions of $\Lambda=\Z[t^{\pm 1}]$.

\begin{example}
If $K$ is the unknot, then $X\cong S^1\x\R^2$ and hence $\tilde X\cong\R^3$.
In this case $\left<A,B\right>=\sum_{n=-\infty}^\infty\lk(A,t^nB)t^n$.
\end{example}

\begin{lemma} \label{cochran-pairing}
The Cochran pairing is 

(a) well-defined;

(b) sesqui-linear in the sense that 

\ \ (b$_1$) $\left<A+A',\,B\right>=\left<A,B\right>+\left<A',B\right>$ and $\left<A,\,B+B'\right>=\left<A,B\right>+\left<A,B'\right>$;

\ \ (b$_2$) $\big\langle q(t)A,\,B\big\rangle=q(t)\left<A,B\right>=\big\langle A,\,q(t^{-1})B\big\rangle$ for any $q(t)\in\Lambda$.

(c) conjugate-symmetric in the sense that $\left<A,B\right>(t)=\left<B,A\right>(t^{-1})$.
\end{lemma}

This is stated without proof in \cite{Co1}*{\S7}.
We include a proof for completeness.

\begin{proof} First we show that $\left<A,B\right>$ does not depend on the choice of $\zeta$.
Let $\zeta'$ be another chain bounded by $f(t)B$.
Since $H_2(\tilde X)=0$, the cycle $\zeta-\zeta'$ bounds a $3$-chain.
Then $(\zeta-\zeta')\cdot t^nB=0$ for each $n$, and the assertion follows.

Next we note that since $\Lambda$ is a PID, the annihilator ideal of $H_1(\tilde X)$ 
is generated by some $\Delta(t)\in\Lambda$.
Let us show that $\left<\cdot,\cdot\right>$, defined using this $\Delta(t)$ as $f(t)$, is sesqui-linear.
Assertion (b$_1$) follows from the linearity of the intersection pairing.
Assertion (b$_2$) reduces to its special case with $q(t)=t^i$ by using (b$_1$).
This case in turn reduces to the case $q(t)=t$ by using induction.
To see that $\left<tA,B\right>=\left<A,\,t^{-1}B\right>$, it suffices to observe that
$t\zeta\cdot t^nB=\zeta\cdot t^{n-1}B$.
To see that $\left<A,\,t^{-1}B\right>=t\left<A,B\right>$, it suffices to note that
$\sum_{n=-\infty}^\infty(\zeta\cdot t^{n-1}B)t^n=
\sum_{k=-\infty}^\infty(\zeta\cdot t^kB)t^{k+1}$.

Next, we show that $\left<A,B\right>$ does not depend on the choice of $f(t)$.
We have $f(t)=\Delta(t)q(t)$ for some $q(t)\in\Lambda$.
If $\Delta(t)A$ bounds a $2$-chain $\xi$ in $\tilde X$, then $f(t)A$ bounds the chain $\zeta:=q(t)\xi$.
Hence $\frac1{f(t)}\sum_{n=-\infty}^\infty(\zeta\cdot t^nB)t^n=
\frac1{q(t)\Delta(t)}\sum_{k=-\infty}^\infty\big(q(t)\xi\cdot t^nB\big)t^n=\frac1{q(t)}\big\langle q(t)A,\,B\big\rangle=\left<A,B\right>$,
where $\left<\cdot,\cdot\right>$ is defined by means of $\Delta(t)$.

Finally, to prove (c) we may use $f(t)=\Delta(t)$.
Then $P:=\Delta(t)A$ bounds a chain $\zeta$, and $Q:=\Delta(t)B$ bounds a chain $\xi$.
Using (b$_2$), it suffices to show that $\left<P,Q\right>(t)=\left<Q,P\right>(t^{-1})$.
It is not hard to see that $\zeta\cdot Q=P\cdot\xi$,%
\footnote{Indeed, $P\x\{0\}$ and $Q\x\{0\}$ bound the chains $\zeta_t:=P\x [0,t]+\zeta\x\{t\}$ 
and $\xi_t:=Q\x [0,t]+\xi\x\{t\}$ in $\tilde X\x I$ for each $t\in I$. 
It is easy to see that $\zeta\cdot Q=\zeta_{1/3}\cdot\xi_{2/3}=\zeta_1\cdot\xi_{2/3}=P\cdot\xi$.}
and we have $P\cdot\xi=(-1)^{1\cdot 2}\xi\cdot P=\xi\cdot P$.
Hence $\zeta\cdot t^nQ=t^n\xi\cdot P=\xi\cdot t^{-n}P$, and it follows that
$\left<P,Q\right>(t)=\left<Q,P\right>(t^{-1})$.
\end{proof}

\begin{corollary} \label{bordism}
Let $A$, $A'$ and $B$ be 1-cycles in $\tilde X$ such that $A$ and $A'$ are disjoint from $t^nB$ for all $n$,
and $A-A'$ bounds a 2-chain $\zeta$ in $\tilde X$ which is disjoint from $t^nB$ for all $n\ne 0$.
Then 

(a) $\left<A,B\right>-\left<A',B\right>=\zeta\cdot B$;

(b) $\left<B,A\right>-\left<B,A'\right>=B\cdot\zeta$.
\end{corollary}

This statement appears (with minor errors) in \cite{PY}*{Lemma 4.5}, \cite{TY}*{Lemma 2.1(i)}.

\begin{proof} Part (a) follows from Lemma \ref{cochran-pairing}(b$_1$).
Part (b) follows from (a) and Lemma \ref{cochran-pairing}(c), using that $B\cdot\zeta=(-1)^{1\cdot 2}\zeta\cdot B=\zeta\cdot B$.
\end{proof}

\begin{remark} \label{pairing-remark}
(a) As observed in \cite{PY}, the image of the Cochran pairing in the quotient $\Q(t)/\Lambda$ is nothing 
but the Blanchfield pairing (see \cite{Hi}). 

(b) Let $Q$ be a knot in $S^3\but K$ such that $\lk(K,Q)=0$, and let $Q_+$ be a parallel pushoff of $P$.
(By {\it parallel} we mean that $\lk(Q,Q_+)=0$.)
Then every lift $\tilde Q$ of $Q$ in $\tilde X$ uniquely determines a lift $\tilde Q_+$ of $Q_+$ in $\tilde X$.
Clearly $\big\langle\tilde Q,\,\tilde Q_+\big\rangle$ does not depend on the choice of the original lift $\tilde Q$.
Hence it is an invariant of the link $(Q,K)$, known as its {\it Kojima's $\eta$-function} \cite{KY}.

(c) Cochran found an expression \cite{Co1}*{Theorem 7.1} for Kojima's $\eta$-function in terms of his derived 
invariants $\beta^i$, which are defined geometrically via Seifert surfaces, without using the infinite cyclic cover.
A similar expression for the Cochran pairing is obtained in \cite{TY}*{Theorem 1.4}, but it is highly asymmetric
with respect to the involution $t\mapsto t^{-1}$ of $\Lambda$.
Below we obtain a more symmetric version (Theorem \ref{cochran expansion}), which generalizes Cochran's original 
formula for Kojima's $\eta$-function.
\end{remark}

\section{Przytycki--Yasuhara theorem} \label{przytycki-yasuhara}

Next let $F$ be a Seifert surface for $K$, and let us fix a lift $\tilde F$ of $F\but K$ in $\tilde X$.
Given a $1$-cycle $Q$ in $S^3\but F$, let $\tilde Q$ denote its lift in $\tilde X$ lying between 
$\tilde F$ and $t\tilde F$.
Given disjoint $1$-cycles $P$, $Q$ in $S^3\but F$, we define $\left<P,Q\right>_F:=\big\langle\tilde P,\tilde Q\big\rangle$.

Let $\phi^\epsilon\:F\to S^3\but F$, where $\epsilon=\pm1$, be the positive and negative parallel pushoffs 
(the sign is determined the co-orientation of $F$).
For a $1$-cycle $Z$ in $F$ let $Z^\epsilon$ denote the $1$-cycle $\phi^\epsilon(Z)$ in $S^3\but F$.
(When $Z$ is a knot, $Z^+$ is not to be confused with the parallel pushoff $Z_+$ of $Z$.
In general, $\lk(Z^+,Z)$ may be nonzero.)
Also, for a homology class $z\in H_1(F)$ let $z^\epsilon$ denote $\phi^\epsilon_*(z)\in H_1(S^3\but F)$.
The following is well-known.

\begin{lemma} \label{seifert-intersection}
Given $\alpha,\beta\in H_1(F)$, we have $\lk(\alpha^+,\beta)-\lk(\beta^+,\alpha)=\alpha\cdot\beta$.
\end{lemma}

\begin{proof}
This follows from $\lk(\beta^+,\alpha)=\lk(\alpha,\beta^+)=\lk(\alpha^-,\beta)$.
\end{proof}

Let $A_1,\dots,A_{2g}\subset F$ be a {\it symplectic basis}, that is, a collection of oriented smoothly embedded 
circles which represent a basis of $H_1(F)$ and are pairwise disjoint from each other, except that each $A_{2i-1}$ 
meets $A_{2i}$ transversely in a single point, with the positive sign (i.e.\ so that $A_{2i-1}\cdot A_{2i}=1$, 
and hence $A_{2i}\cdot A_{2i-1}=-1$).
We may then view $F$ as a disk $D$ with bands, where the core of the $i$th band is an arc lying in $A_i$
and the rest of $A_i$ lies in $D$.
Let $D_1,\dots,D_{2g}\subset S^3$ be a collection smooth $2$-disks such that each $B_i:=\partial D_i$ lies 
in $S^3\but F$, each $D_i$ intersects $F$ transversely in an arc disjoint from $D$ and meeting $A_i$ in a single point, 
and $D_i\cap D_j=\emptyset$ for $i\ne j$.
We orient the $B_i$ so that $\lk(A_i,B_j)=\delta_{ij}$.
We will refer to the knots $B_1,\dots,B_{2g}\subset S^3\but F$ as the {\it meridians} of the symplectic basis
$A_1,\dots,A_{2g}$.

Let us note that $[A_i]=\theta([B_i])$, where 
$\theta\:H_1(S^3\but F)\xr{\simeq}H^1(F)\xr\simeq\Hom\big(H_1(F),\Z\big)\xr{\simeq} H_1(F)$
is the composition of the Alexander duality and the isomorphisms given by the universal coefficient formula and by
the non-degeneracy of the intersection pairing.
In particular, since the $[A_i]$ form a basis of $H_1(F)$, the $[B_i]$ form a basis of $H_1(S^3\but F)$.
Moreover, for any $1$-cycle $Q$ in $S^3\but F$ the coefficients of its expression in this basis, 
$[Q]=\sum_{i=1}^{2g}c_i[B_j]$, are given by the following observation: $\lk(Q,A_j)=\sum_{i=1}^{2g}c_i\lk(B_i,A_j)=c_j$.
In particular, we have $[A_i^+]=\sum_{j=1}^{2g} c_{ij}[B_j]$, where $c_{ij}=\lk(A_i^+,A_j)$.
The matrix $V:=\big(\lk\big(A_i^+,A_j)\big)$ is called the Seifert matrix of $F$ with respect to the chosen basis
of $H_1(F)$.
Lemma \ref{seifert-intersection} immediately yields the following well-known fact:

\begin{corollary} \label{seifert-intersection2}
$V-V^T$ is the intersection matrix of $H_1(F)$ with respect to the chosen basis, that is, the block diagonal matrix 
where all diagonal blocks are $\begin{pmatrix} 0 & 1\\ -1 & 0\end{pmatrix}$.
\end{corollary}

Corollary \ref{seifert-intersection2} implies that $\det(V-V^T)=1$.
In particular, the matrix $tV-V^T$ is nonsingular (in other words, invertible over the quotient field $\Q(t)$ of $\Z[t]$).%
\footnote{In fact, $\det(tV-V^T)$ is the Alexander polyniomial of $K$ (see Example \ref{kauffman}).}

\begin{proposition}[Przytycki--Yasuhara] \label{PY1} \

(a) \cite{PY}*{Lemma 4.6} $\Big(\left<B_i,B_{j+}\right>_F\Big)_{1\le i,j\le 2g}=(1-t)(tV-V^T)^{-1}$.
\smallskip

(b) $\Big(\left<B_i,B_{j+}\right>_F\Big)_{1\le i,j\le 2g}=(x-x^{-1})(-xV+x^{-1}V^T)^{-1}$, where $t=x^2$.
\end{proposition}

Beware that $\left<B_i,B_{j+}\right>_F$ is written as $\left<B_i,B_j\right>_F$ (in a different notation) 
in \cite{PY}, which makes no sense when $i=j$. 
The proof of (a) in \cite{PY} also lacks clarity in treating pushoffs of $B_j$'s.
For completeness we include a more accurate exposition of the proof of (a).

\begin{proof}[Proof. (a)] Let $N$ be a regular neighborhood of $F$ in $S^3$
and let $N^+$ be a regular neighborhood of $N$ in $S^3$.
We may assume that each $A_i^+$ and $A_i^-$ lie in $S^4\but N^+$, each $B_i$ lies in $N^+\but N$ and 
each $B_{i+}$ lies in $N$ and moreover bounds a disk in $N$.
Let $\tilde B_i$, $\tilde B_{i+}$ and $\tilde A_i^\pm$ be the lifts of $B_i$, $B_{i+}$
and $A_i^\pm$ in $\tilde X$ lying in the region bounded by $\tilde F$ and $t\tilde F$.

As noted above, each $[A_i^+]=\sum_{j=1}^{2g}\lk(A_i^+,A_j)[B_j]$ in $H_1(S^4\but F)$, and hence also in 
$H_1(S^3\but N)$; similarly, each $[A_i^-]=\sum_{j=1}^{2g}\lk(A_i^-,A_j)[B_j]$.
Let $\theta_i^\pm$ be a $2$-chain in $S^3\but N$ bounded by $A_i^\pm-\sum_{j=1}^{2g}\lk(A_i^\pm,A_j)B_j$.
Let $\tilde\theta_i^\pm$ be its lift in $\tilde X$ lying in the region bounded by $\tilde F$ and $t\tilde F$.
Thus $\partial\tilde\theta_i^\pm=\tilde A_i^\pm-\sum_{j=1}^{2g}\lk(A_i^\pm,A_j)\tilde B_j$.
Let $\eta^\pm_i$ be a $2$-chain in $X=S^3\but K$ bounded by $A_i^\pm-A_i$, with support in an annulus 
meeting $F$ only in $A_i$.
Let $\tilde\eta^\pm_i$ be the lift of $\eta^\pm_i$ in $\tilde X$ that intersects $t\tilde F$,
and let $\tilde\eta_i=\tilde\eta^+_i-\tilde\eta^-_i$.
Thus $\partial\eta_i=t\tilde A_i^+-\tilde A_i^-$.

Let $\xi_i=\tilde\eta_i-(t\tilde\theta_i^+-\tilde\theta_i^-)$.
Then $\partial\xi_i=\sum_{j=1}^{2g}\big(\lk(A_i^+,A_j)t-\lk(A_i^-,A_j)\big)\tilde B_j$.
Since $\lk\big(A_i^-,A_j)=\lk\big(A_i,A_j^+)=\lk(A_j^+,A_i)$, we obtain that
\[\begin{pmatrix}\partial\xi_1\\ \vdots\\ \partial\xi_{2g}\end{pmatrix}=
(tV-V^T)\begin{pmatrix}\tilde B_1\\ \vdots\\ \tilde B_{2g}\end{pmatrix}.\]

Let $M=tV-V^T$.
As noted above, $M$ is invertible over the quotient field $\Q(t)$ of $\Z[t]$.
Moreover, by Lemma \ref{adjugate} $(\det M)M^{-1}$ is the adjugate matrix $\adj M$, whose entries belong to $\Z[t]$.
Let us multiply both sides of the previous equation by $\adj M$ on the left:
\[\adj(M)\begin{pmatrix}\partial\xi_1\\ \vdots\\ \partial\xi_{2g}\end{pmatrix}=
\det(M)\begin{pmatrix}\tilde B_1\\ \vdots\\ \tilde B_{2g}\end{pmatrix}.\]
Writing $\mu_{ij}$ for the entries of $\adj M$, we get that $\det(M)\tilde B_i=\partial\zeta_i$,
where $\zeta_i=\sum_{j=1}^{2g}\mu_{ij}\xi_j$.
Therefore
\[\left<B_i,B_{j+}\right>_F=\frac1{\det M}\sum_{n=-\infty}^\infty(\zeta_i\cdot t^n\tilde B_{j+})t^n.\]

In order to find $\zeta_i\cdot t^nB_{j+}$, let us first compute $\xi_i\cdot t^nB_{j+}$.
Since each $\theta_i^\pm$ lies in $S^3\but N$ and each $B_{j+}$ lies in $N$, we have 
$\tilde\theta_i^\pm\cdot t^n\tilde B_{j+}=0$ for all $n\in\Z$.
Also it is easy to see that $\tilde\eta_i\cdot t^n\tilde B_{j+}=0$ for $n\ne 0,1$.
Since $B_{j+}$ bounds a disk in $N$ and $A_i^\pm$ lies in $S^3\but N$, we have $\lk(A_i^\pm,B_{j+})=0$.
Hence $\eta^\pm_i\cdot B_{j+}=\lk(A_i^\pm,B_{j+})-\lk(A_i,B_{j+})=-\lk(A_i,B_{j+})=-\delta_{ij}$
(the Kronecker symbol).
Then $\tilde\eta_i\cdot\tilde B_{j+}=-\tilde\eta^-_i\cdot\tilde B_{j+}=-\eta^-_i\cdot B_{j+}=\delta_{ij}$ and
$\tilde\eta_i\cdot t\tilde B_{j+}=\tilde\eta^+_i\cdot t\tilde B_{j+}=\eta^+_i\cdot B_{j+}=-\delta_{ij}$.
Thus \[\xi_i\cdot t^n\tilde B_{j+}=\tilde\eta_i\cdot t^n\tilde B_{j+}=\begin{cases}
\delta_{ij}&\text{if }n=0\\
-\delta_{ij}&\text{if }n=1\\
0&\text{if }n\ne 0,1.
\end{cases}\]
Hence \[\zeta_i\cdot t^n\tilde B_{j+}=\sum_{k=1}^{2g}\mu_{ik}\xi_k\cdot t^n\tilde B_{j+}=
\mu_{ij}\xi_j\cdot t^n\tilde B_{j+}=\begin{cases}
\mu_{ij}&\text{if }n=0\\
-\mu_{ij}&\text{if }n=1\\
0&\text{if }n\ne 0,1.
\end{cases}\]
Thus $\left<B_i,B_{j+}\right>_F=(1-t)\dfrac{\mu_{ij}}{\det M}$, where $\dfrac{\mu_{ij}}{\det M}$ is the $(i,j)$-entry of
$\dfrac{\adj M}{\det M}=M^{-1}$.
\end{proof}

\begin{proof}[(b)] We have
$(1-x^2)(x^2V-V^T)^{-1}=(x^{-1}-x)x\cdot (-x^{-1})(-xV+x^{-1}V^T)^{-1}$,
so the assertion follows from (a). 
\end{proof}

The following is well-known:

\begin{lemma} \label{adjugate} 
For a square matrix $M$ over a ring, $M\adj M=\det(M)I$, where the adjugate matrix $\adj M$
is the transpose of the cofactor matrix $\Big((-1)^{i+j}\det(M_{ij})\Big)_{ij}$, each $M_{ij}$ being the
matrix obtained from $M$ by deleting the $i$th row and the $j$th column.
\end{lemma}

\begin{theorem}[Przytycki--Yasuhara] \label{PY2}
Let $P$ and $Q$ be disjoint smooth knots in $S^3\but F$, and let us write 
$L_P=\big(\lk(P,A_1),\dots,\lk(P,A_{2g})\big)$. 
Then

(a) \cite{PY}*{Theorem 4.1} $\left<P,Q\right>_F=\lk(P,Q)+(1-t)L_P(tV-V^T)^{-1}L_Q^T$.
\smallskip

(b) $\left<P,Q\right>_F=\lk(P,Q)+(x-x^{-1})L_P(-xV+x^{-1}V^T)^{-1}L_Q^T$, where $t=x^2$.
\end{theorem}

The statement of (a) is applied in \cite{TY} to obtain a remarkable factorization theorem, 
which is further applied in \cite{M21}.

The proof of (a) in \cite{PY} lacks clarity in treating pushoffs of $B_j$'s.
For completeness we include a more accurate exposition of the proof of (a).

\begin{proof}[Proof. (a)] Let $N$ be a regular neighborhood of $F$ in $S^3$.
We may assume that $P$, $Q$ and the $B_i$ lie in $S^3\but N$, but the parallel pushoffs $B_{i+}$ lie in $N$
and moreover bound disks in $N$.
Since $[P]=\sum_{i=1}^{2g}\lk(P,A_i)[B_i]$ in $H_1(S^3\but F)$, and hence also in $H_1(S^3\but N)$,
there is a $2$-chain $\zeta$ in $S^3\but N$ such that $\partial\zeta=P-\sum_{i=1}^{2g}\lk(P,A_i)B_i$.
Similarly there is a $2$-chain $\xi$ in $S^3\but F$ such that $\partial\xi=Q-\sum_{j=1}^{2g}\lk(Q,A_j)B_{j+}$.
Let $\tilde P$, $\tilde Q$, $\tilde B_i$, $\tilde B_{j+}$, $\tilde\zeta$ and $\tilde\xi$ be the lifts of 
$P$, $Q$, $B_i$, $B_{j+}$, $\zeta$ and $\xi$ in $\tilde X$ lying in the region bounded by $\tilde F$ and $t\tilde F$.

Since each $B_{j+}$ bounds a disk in $N$ but $P$ lies in $S^3\but N$, we have $\lk(P,B_{j+})=0$.
Hence $\tilde P\cdot\tilde\xi=P\cdot\xi=\lk(P,Q)-\lk\big(P,\sum_{i=1}^{2g}\lk(Q,A_j)B_{j+}\big)=\lk(P,Q)$.
Then by Corollary \ref{bordism}(b) and Lemma \ref{cochran-pairing}(b$_1$)
\[\big\langle\tilde P,\tilde Q\big\rangle=\tilde P\cdot\tilde\xi+\Big\langle\tilde P,\sum_{j=1}^{2g}\lk(Q,A_j)\tilde B_{j+}\Big\rangle=
\lk(P,Q)+\sum_{j=1}^{2g}\lk(Q,A_j)\big\langle\tilde P,\tilde B_{j+}\big\rangle.\]
Since $B_{j+}$ lies in $N$ but $\zeta$ lies in $S^3\but N$, we have 
$\tilde\zeta\cdot\tilde B_{j+}=\zeta\cdot B_{j+}=0$.
Hence by Corollary \ref{bordism}(a) and Lemma \ref{cochran-pairing}(b$_1$)
\[\big\langle\tilde P,\tilde B_{j+}\big\rangle=\tilde\zeta\cdot\tilde B_{j+}+\Big\langle\sum_{i=1}^{2g}\lk(P,A_i)\tilde B_i,\tilde B_{j+}\Big\rangle=
\sum_{i=1}^{2g}\lk(P,A_i)\big\langle\tilde B_i,\tilde B_{j+}\big\rangle.\]
Thus \[\big\langle\tilde P,\tilde Q\big\rangle=
\lk(P,Q)+\sum_{i=1}^{2g}\sum_{j=1}^{2g}\lk(P,A_i)\big\langle\tilde B_i,\tilde B_{j+}\big\rangle\lk(A_j,Q).\]
Now the desired assertion follows from Proposition \ref{PY1}(a).
\end{proof}

\begin{proof}[(b)] This follows from (a) similarly to the proof of Proposition \ref{PY1}(b).
\end{proof}

\section{Generalized Cochran invariants}\label{generalized cochran}

Let $K$ be a smooth knot in $S^3$ and $F$ a Seifert surface for $K$.
Let $Z$ be a $1$-cycle in $S^3\but F$.
Since $Z$ is disjoint from $F$, we have $\lk(Z,K)=0$, and consequently $Z$ bounds a $2$-chain $\zeta$ in $S^3\but K$.
We may assume that it meets $F$ transversely, which results in a $1$-cycle $D(Z):=\zeta\cap F$.

\begin{lemma} \label{intersection} \cite{TY}*{1.3} 
(a) $[D(Z)]\in H_1(F)$ is uniquely determined by $[Z]\in H_1(S^3\but F)$.

(b) For a $1$-cycle $W$ in $F$ we have $\lk(Z,W)=D(Z)\cdot W$.
\end{lemma}

\begin{proof}[Proof. (a)] It is easy to see that the assignment $[Z]\mapsto[D(Z)]$ is nothing but a standard geometric
description of the composition $H_1(S^3\but F)\xr{\simeq} H^1(F)\xr{\simeq} H_1(F,\partial F)\xr{\simeq} H_1(F)$ of 
the Alexander duality, the Poincar\'e duality and the inverse of the homomorphism $H_1(F)\to H_1(F,\partial F)$ 
(which is easily seen to be an isomorphism).
\end{proof}

\begin{proof}[(b)]
$\lk(Z,W)=\zeta\cdot W=(\zeta\cap F)\cdot W=D(Z)\cdot W$.
\end{proof}

\begin{proof}[Alternative proof of (a)]
$[D(Z)]=\sum_{i=1}^{2g} c_i\alpha_i$, where $\alpha_1,\dots,\alpha_{2g}$ is 
a symplectic basis for $H_1(F)$.
By applying $\cdot\alpha_{2j}$ and $\cdot\alpha_{2j+1}$ to both sides of this equation 
and using (b) we find that $c_{2j-1}=[D(Z)]\cdot\alpha_{2j}=\lk([Z],\alpha_{2j})$ and 
$c_{2j}=-[D(Z)]\cdot\alpha_{2j-1}=-\lk([Z],\alpha_{2j-1})$.
\end{proof}

For a $1$-cycle $Z$ in $S^3\but F$ let $Z^{(0)}=Z$ and $Z^{(n)}=D(Z^{(n-1)+})$ for $n\ge 1$.
Also, we set $Z^\pm=Z$. Given disjoint $1$-cycles $P$, $Q$ in $S^3\but F$, let 
$\beta^{k,l}_F(P,Q)=\lk\big(P^{(k)+},\,Q^{(l)}\big)$.

Clearly, $\beta_F^{0,0}(P,Q)=\lk(P,Q)$.

\begin{lemma} $\beta^{k,l}_F(P,Q)$ is well-defined, and for $(k,l)\ne (0,0)$ it depends only on 
$[P],[Q]\in H_1(S^3\but F)$.
\end{lemma}

\begin{proof} First we note that Lemma \ref{intersection}(a) along with the maps 
$\phi^\eps_*\:H_1(F)\to H_1(S^3\but F)$ yield that
$[Q]\in H_1(S^3\but F)$ uniquely determines $[Q^{(n)}]\in H_1(F)$ for each $n\ge 1$ and 
$[Q^{(n)+}],\,[Q^{(n)-}]\in H_1(S^3\but F)$ for each $n\ge 0$.

Now if $l>0$, then $[P^{(l)}]\in H_1(F)$ is defined, and has a well-defined linking number with 
$[Q^{(k)+}]\in H_1(S^3\but F)$.
If $k>0$, then $[Q^{(k)}]\in H_1(F)$ is defined, and has a well-defined linking number with 
$[P^{(l)-}]\in H_1(S^3\but F)$.
But clearly $\beta^{k,l}_F(P,Q)=\lk\big(P^{(k)},\,Q^{(l)-}\big)$.
\end{proof}

\begin{remark} 
(a) $\beta_F^{n,n}(Q,Q)$ for $n\ge 1$ is nothing but Cochran's derived invariant $\beta^n(Q,K)$ \cite{Co1}.

(b) $\beta_F^{n,0}(P,Q)$ coincides with the namesake invariant of \cite{TY}.

(c) $\beta_F^{k,l}(Q,Q)$ for $k,l\ge 1$ were studied in \cite{GL} (see also \cite{Co1}*{\S8}).
In particular, they were shown in \cite{GL}*{\S\S4,5} to be determined by Cochran's derived invariants 
$\beta^n(Q,K)=\beta_F^{n,n}(Q,Q)$, $n\ge 1$.
\end{remark}

\begin{lemma} \label{intersection2} 
(a) If $k\ge 0$ and $l\ge 1$, then $\beta^{k,l}_F(P,Q)=P^{(k+1)}\cdot Q^{(l)}$.
\smallskip

(b) If $k\ge 1$, then $\beta^{k,0}_F(P,Q)=-P^{(k)}\cdot Q^{(1)}$.
\end{lemma}

\begin{proof}[Proof. (a)]
By Lemma \ref{intersection}(b) $\beta^{k,l}_F(P,Q)=D(P^{(k)+})\cdot Q^{(l)}=P^{(k+1)}\cdot Q^{(l)}$.
\end{proof}

\begin{proof}[(b)] This follows from (a) and the obvious relation $\beta^{k,0}(P,Q)=\beta^{0,k}(Q,P)$.
\end{proof}

\begin{lemma} \label{pascal} Suppose that $k,l\ge 0$.
\smallskip

(a) $\beta^{k,\,l+1}_F(P,Q)=-\beta^{l,\,k+1}_F(Q,P)$.
\smallskip

(b) $\beta^{k+1,\,0}_F(P,Q)=-\beta^{k,1}_F(P,Q)$.
\smallskip

(c) $\beta^{k,\,k+1}_F(Q,Q)=0$.
\smallskip

(d) If $l\ge 1$, then $\beta^{k,l}_F(P,Q)-\beta^{l,k}_F(Q,P)=\beta^{k-1,\,l}_F(P,Q)$.
\smallskip

(e) If $l\ge 1$, then $\beta^{k+1,\,l}_F(P,Q)+\beta^{k,\,l+1}_F(P,Q)=\beta^{k,l}_F(P,Q)$.
\end{lemma}

Parts (e) and (b) imply that all the $\beta^{k,l}_F(P,Q)$ are determined by $\beta^{n,0}_F(P,Q)$;
and also by $\beta^{n,n}_F(P,Q)$ and $\beta^{n,\,n+1}_F(P,Q)$ (compare \cite{GL}*{\S5}).

\begin{proof}[Proof. (a)] This follows from Lemma \ref{intersection2}(a).
\end{proof}

\begin{proof}[(b)] 
This follows by comparing parts (a) and (b) of Lemma \ref{intersection2}.
\end{proof}

\begin{proof}[(c)]
This follows from (a). 
Alternatively, since $D(Q)$ lies on a $2$-chain bounded by $Q$, we have $\lk\big(D(Q),\,Q\big)=0$.
\end{proof}

\begin{proof}[(d)]
Given $1$-cycles $Z$, $W$ in $F$, by Lemma \ref{seifert-intersection} we have $\lk(Z^+,W)-\lk(W^+,Z)=Z\cdot W$.
It remains to apply this formula to $Z=P^{(k)}$ and $W=Q^{(l)}$ and use Lemma \ref{intersection2}(a).
\end{proof}

\begin{proof}[(e)] By (d) $\beta^{k+1,\,l}_F(P,Q)+\beta^{l,\,k+1}_F(Q,P)=\beta^{k,l}_F(P,Q)$, where
$\beta^{l,\,k+1}_F(Q,P)=-\beta^{k,\,l+1}_F(P,Q)$ by (a).
\end{proof}

\begin{theorem} \label{cochran expansion} Let $P$ and $Q$ be disjoint $1$-cycles in $S^3\but F$.

(a) Then $\displaystyle\left<P,Q\right>_F=\sum_{n=0}^\infty\Big(\beta_F^{n,n}(P,Q)-(1-t)\beta_F^{n,\,n+1}(P,Q)\Big)(1-t)^n(1-t^{-1})^n$.

(b) Let $x^2=t$. Then
\[\left<P,Q\right>_F=\beta_F^{0,0}(P,Q)+
\sum_{n=1}^\infty(-1)^n\Big(x\beta_F^{n,n}(Q,P)-x^{-1}\beta_F^{n,n}(P,Q)\Big)(x-x^{-1})^{2n-1}.\]
\end{theorem}

In connection with (b), let us note that $\beta_F^{n,n}(Q,P)=-\beta_F^{n-1,\,n+1}(P,Q)$ for $n\ge 1$ by Lemma \ref{pascal}(a).

Each of the two assertions of Theorem \ref{cochran expansion} is a generalization of \cite{Co1}*{Theorem 7.1} (concerning (a),
see Lemma \ref{pascal}(c)).
It is shown in \cite{TY}*{Theorem 1.4} that \[\left<P,Q\right>_F=\sum_{n=0}^\infty(-1)^n\beta_F^{n,0}(P,Q)(t-1)^n.\]

\begin{proof}[Proof. (a)] Let $X=S^3\but K$, let $\tilde X$ be its infinite cyclic cover and $t$ be the positive generator
of the group of covering transformations of $\tilde X$.
Let $F_1$ be an auxiliary positive pushoff of $F$.
(The reason why we do not simply use $F$ instead of $F_1$ will be clear closer to the end of the proof.)
Let $F_1^+$ be a positive pushoff of $F_1$ and let $\tilde F_1$ and $\tilde F_1^+$ be the lifts of $F_1\but K$ and $F_1^+\but K$
respectively in $\tilde X$ that lie near $\tilde F$ (so, between $\tilde F$ and $t\tilde F$).
We may assume that $\tilde P$ and $\tilde Q$ lie between $\tilde F_1^+$ and $t\tilde F$.
Let $\zeta$ and $\xi$ be some $2$-chains bounded by $P$ and $Q$ in $S^3\but K$.
Let $P_1^+=\zeta\cap F_1^+$ and $Q_1=\xi\cap F_1$, which we assume to be transverse intersections, and let us consider 
their lifts $\tilde P_1^+\subset\tilde F_1$ and $\tilde Q_1\subset\tilde F_1$.
Let $\tilde\zeta_+$ be the $2$-chain in $\tilde X$ lying between $\tilde F_1^+$ and $t\tilde F_1^+$ and projecting onto $\zeta$,
and let $\tilde\xi$ be the $2$-chain in $\tilde X$ lying between $\tilde F_1$ and $t\tilde F_1$ and projecting onto $\xi$.
Then $\partial\tilde\zeta_+=\tilde P-(1-t)\tilde P_1^+$ and $\partial\tilde\xi=\tilde Q-(1-t)\tilde Q_1$.

We have $\tilde P\cdot\tilde\xi=P\cdot\xi=\lk(P,Q)$.
Then by Corollary \ref{bordism}(b) and Lemma \ref{cochran-pairing}(b$_2$) 
\[\big\langle\tilde P,\tilde Q\big\rangle=\big\langle\tilde P,\,\partial\tilde\xi+(1-t)\tilde Q_1\big\rangle=
\lk(P,Q)+(1-t^{-1})\big\langle\tilde P,\tilde Q_1\big\rangle.\]
Next, $\tilde\zeta_+\cdot t\tilde Q_1=\zeta\cdot Q_1=\lk(P,Q_1)$.
Then by Corollary \ref{bordism}(a) and Lemma \ref{cochran-pairing}(b$_2$) 
\[\big\langle\tilde P,\tilde Q_1\big\rangle=t\big\langle\partial\tilde\zeta_++(1-t)\tilde P_1^+,\,t\tilde Q_1\big\rangle=
t\lk(P,Q_1)+(1-t)\big\langle\tilde P_1^+,\tilde Q_1\big\rangle.\]
Thus \[\big\langle\tilde P,\tilde Q\big\rangle=\lk(P,Q)-(1-t)\lk(P,Q_1)+(1-t)(1-t^{-1})\big\langle\tilde P_1^+,\tilde Q_1\big\rangle.\]

Next let $F_2$ be a positive pushoff of $F$ and $F_2^+$ a positive pushoff of $F_2$, both chosen so as to lie
between $F$ and $F_1$.
Let $\eta$ and $\theta$ be some $2$-chains bounded by $P_1^+$ and $Q_1$ in $S^3\but K$.
Let $P_2^+=\eta\cap F_2^+$ and $Q_2=\theta\cap F_2$, which we assume to be transverse intersections, and let us consider 
their lifts $\tilde P_2^+\subset\tilde F_2^+$ and $\tilde Q_2\subset\tilde F_2$.
Then similarly to the above 
\[\big\langle\tilde P_1^+,\tilde Q_1\big\rangle=
\lk(P_1^+,Q_1)-(1-t)\lk(P_1^+,Q_2)+(1-t)(1-t^{-1})\big\langle\tilde P_2^+,\tilde Q_2\big\rangle.\]

It is easy to see that $\lk(P,Q_1)=\lk(P,Q^{(1)})$, $\lk(P_1^+,Q_1)=\lk(P^{(1)+},Q^{(1)})$ and, crucially,%
\footnote{If we did the same procedure in the opposite order, that is, first applied
$\tilde\zeta_+\cdot\tilde Q=\zeta\cdot Q=\lk(P,Q)$
and then 
$\tilde P_1^+\cdot\tilde\xi=P_1^+\cdot\xi=\lk(P_1,Q)$,
we would similarly obtain
$\big\langle\tilde P,\tilde Q\big\rangle=\lk(P,Q)+(1-t)\lk(P_1^+,Q)+(1-t)(1-t^{-1})\big\langle\tilde P_1^+,\tilde Q_1\big\rangle$
and $\big\langle\tilde P_1^+,\tilde Q_1\big\rangle=
\lk(P_1^+,Q_1)+(1-t)\lk(P_2^+,Q_1)+(1-t)(1-t^{-1})\big\langle\tilde P_2^+,\tilde Q_2\big\rangle$.
But it is not true that $\lk(P_2^+,Q_1)=\lk(P^{(2)+},Q^{(1)})$.}
$\lk(P_1^+,Q_2)=\lk(P^{(1)+},Q^{(2)})$.
Then by continuing in the same fashion we obtain the desired formula.
\end{proof}

\begin{proof}[(b)] It is easy to see that $(1-t)(1-t^{-1})=(1-t)+(1-t^{-1})$.
On the other hand, $\beta_F^{n+1,\,n+1}(P,Q)-\beta_F^{n,\,n+1}(P,Q)=\beta_F^{n+1,\,n+1}(Q,P)$ by Lemma \ref{pascal}(d).
Using these facts, part (a) implies
\[\left<P,Q\right>_F=\beta_F^{0,0}(P,Q)+
\sum_{n=0}^\infty\Big((1-t)\beta_F^{n+1,\,n+1}(Q,P)+(1-t^{-1})\beta_F^{n+1,\,n+1}(P,Q)\Big)(1-t)^n(1-t^{-1})^n.\]
Also $1-t=1-x^2=-x(x-x^{-1})$ and similarly $1-t^{-1}=x^{-1}(x-x^{-1})$, whence $(1-t)(1-t^{-1})=-(x-x^{-1})^2$.
This yields the desired formula.
\end{proof}

Let us consider the isomorphism $h\:H_1(S^3\but F)\to H_1(F)$ from the first proof of Lemma~\ref{intersection}(a) 
and the composition $f\:H_1(F)\xr{\phi^+_*}H_1(S^3\but F)\xr{h}H_1(F)$.
Clearly, for $n\ge 1$ we have $[Q^{(n)}]=f([Q^{(n-1)}])$, and consequently $[Q^{(n)}]=f^{n-1}(q)$, where $q=h([Q])$.
Thus by Lemma \ref{intersection2} for $l\ge 1$ we have $\beta^{k,l}_F(P,Q)=f^k(p)\cdot f^{l-1}(q)$, where $p=h([P])$. 

Here is a more explicit version of this computation (compare \cite{Co1}*{\S8}, \cite{TY}*{\S4}).

\begin{proposition} \label{cochran-computation}
Let $b=(\alpha_1\dots,\alpha_{2g})$ be a symplectic basis for $H_1(F)$.
Let $V=\big(\lk(\alpha_i^+,\alpha_j)\big)$ be its Seifert matrix and $J=(\alpha_i\cdot\alpha_j)$ 
be its intersection matrix.
For a $1$-cycle $Z$ in $S^3\but F$ let $L_Z=\big(\lk([Z],\alpha_1),\dots,\lk([Z],\alpha_{2g})\big)$. 
\smallskip

(a) If $k\ge 0$ and $l\ge 1$, then $\beta^{k,l}_F(P,Q)=L_P(J^TV)^kJ(V^TJ)^{l-1}L_Q^T$.
\smallskip

(b) If $k\ge 1$, then $\beta^{k,\,0}_F(P,Q)=-L_P(J^TV)^{k-1}JL_Q^T$.
\end{proposition}

Since $V-V^T=J=-J^T$, we have $JV^T=J(V+J^T)=JV+I=I-J^TV$.
Using this, it is easy to get an explicit expression for each $\beta^{k,l}$ in terms of the $\beta^{k,0}$.

\begin{proof} 
Let $Z$ be a $1$-cycle in $S^3\but F$.
By the second proof of Lemma \ref{intersection}
$[D(Z)]=bJL_Z^T\in H_1(F)$.
This can also be written as $[D(Z)]=L_ZJ^Tb^T$.
Let $b_+=(\alpha_1^+,\dots,\alpha_{2g}^+)$.
Then $[D(Z)^+]=L_ZJ^Tb_+^T\in H_1(S^3\but F)$.
It follows that $L_{D(Z)^+}=L_ZJ^TV$.
Therefore $L_{Z^{(n)}}=L_{Z^{(n-1)}}J^TV=L_Z(J^TV)^n$.
Hence $[Z^{(n)}]=[D(Z^{(n-1)})^+]=L_{Z^{(n-1)}}J^Tb_+^T=L_Z(J^TV)^{n-1}J^Tb^T\in H_1(S^3\but F)$.

For $k\ge 0$ and $l\ge 1$ by Lemma \ref{intersection2}(a)
$\beta_F^{k,l}(P,Q)=P^{(k+1)}\cdot Q^{(l)}$.
Therefore $\beta_F^{k,l}(P,Q)=\big(L_P(J^TV)^kJ^T\big)J\big(L_Q(J^TV)^{l-1}J^T\big)^T$.
Since $J^TJ=I$, we get the assertion of (a).

For $k\ge 1$ by Lemma \ref{pascal}(b) $\beta^{k,\,0}_F(P,Q)=-\beta^{k-1,1}_F(P,Q)=-L_P(J^TV)^{k-1}JL_Q^T$.
\end{proof}

\begin{proposition} Let $\alpha_1\dots,\alpha_{2g}$ be a symplectic basis for $H_1(F)$.
Let $V=\big(\alpha_i^+,\alpha_j\big)$ be its Seifert matrix and $J=(\alpha_i\cdot\alpha_j)$ 
be its intersection matrix.
Then
\[(-xV+x^{-1}V^T)^{-1}=\sum_{n=1}^\infty(-1)^n\Big((J^TV)^{n-1}J^T(xV^T-x^{-1}V)J(V^TJ)^{n-1}\Big)(x-x^{-1})^{2n-2}.\]
\end{proposition}

Compare \cite{PY}*{Corollary 4.3}.

\begin{proof} Let $P$ and $Q$ be disjoint $1$-cycles in $S^3\but F$.
Then by Theorem \ref{PY2}(b) we have
$\left<P,Q\right>_F-\lk(P,Q)=(x-x^{-1})L_P(-xV+x^{-1}V^T)^{-1}L_Q^T$.
By Theorem \ref{cochran expansion}(b)
\[\left<P,Q\right>_F-\lk(P,Q)=
\sum_{n=1}^\infty(-1)^n\Big(x\beta_F^{n,n}(Q,P)-x^{-1}\beta_F^{n,n}(P,Q)\Big)(x-x^{-1})^{2n-1}.\]
For $n\ge 1$ by Lemma \ref{pascal}(a) $\beta_F^{n,n}(Q,P)=-\beta_F^{n-1,\,n+1}(P,Q)$.
Then by Proposition \ref{cochran-computation}(a)
$x\beta_F^{n,n}(Q,P)-x^{-1}\beta_F^{n,n}(P,Q)=L_P(J^TV)^{n-1}(-xJV^TJ-x^{-1}J^TVJ)(V^TJ)^{n-1}L_Q^T$.
Since $J^T=-J$, this simplifies as
\[x\beta_F^{n,n}(Q,P)-x^{-1}\beta_F^{n,n}(P,Q)=L_P(J^TV)^{n-1}J^T(xV^T-x^{-1}V)J(V^TJ)^{n-1}L_Q^T.\]
Thus $L_PM_0L_Q^T=\sum_{n=1}^\infty L_PM_nL_Q^T$, where $M_0=(-xV+x^{-1}V^T)^{-1}$
and \[M_n=(-1)^n\Big((J^TV)^{n-1}J^T(xV^T-x^{-1}V)J(V^TJ)^{n-1}\Big)(x-x^{-1})^{2n-2}\]
for $n\ge 1$.
Let $[B_1],\dots,[B_{2g}]\in H_1(S^3\but F)$ be the ``dual'' basis, satisfying $\lk(\alpha_i,[B_j])=\delta_{ij}$
(see \S\ref{cochran-pairing-section}).
Then $L_{B_j}M_kL_{B_i}^T$ is the $(i,j)$-entry of $M_k$, for each $k$.
Hence we get that $\sum_{n=1}^\infty M_n$ converges and has the same entries as $M_0$.
\end{proof}

\section{Conway potential function} \label{cpf}

Let $L=(L_1,\dots,L_n)$ be a colored smooth link in $S^3$.
A {\it C-complex} for $L$ is a collection $F=(F_1,\dots,F_n)$ of smoothly embedded oriented surfaces in $S^3$
such that
\begin{itemize}
\item each $F_i$ is connected and each $\partial F_i=L_i$ (as oriented $1$-manifolds);
\item if $i\ne j$, then $F_i$ and $F_j$ intersect transversely, and $F_i\cap F_j$ is 
a disjoint union of {\it clasp arcs}, that is, arcs $J_{ijk}$ such that each $J_{ijk}$ 
has one boundary point in $L_i$ and another in $L_j$, and is otherwise disjoint from $L$;
\item $F_i\cap F_j\cap F_k=\emptyset$ for every pairwise distinct $i,j,k$;
\item $|F|:=F_1\cup\dots\cup F_n$ is connected.
\end{itemize}

\begin{theorem} \label{cimasoni0} \cite{Coo1}, \cite{Ci0} Every colored smooth link has a C-complex.
\end{theorem}

See also Lemmas \ref{C-complex} and \ref{C-complex2} for a proof for $2$-component links with $\lk=0,1$.

\begin{example}
It is well-known that the Borromean rings $(L_1,L_2,L_3)$ do not bound Seifert surfaces 
$(F_1,F_2,F_3)$ such that $F_1\cap F_2\cap F_3=\emptyset$ and each $L_i\cap F_j=\emptyset$ for $i\ne j$.
However, it is not hard to see that they bound embedded disks $(D_1,D_2,D_3)$ such that 
$D_1\cap D_3=\emptyset$ and $D_2$ meets each of $D_1$ and $D_3$ in two clasp arcs.
\end{example}

A clasp arc $J_{ijk}$ is {\it positive} if its orientation (that is, the orientation of the intersection%
\footnote{Let us note that this orientation is reversed if $F_i$ and $F_j$ are interchanged.}
between the oriented surfaces $F_i$ and $F_j$ in the oriented $3$-manifold $S^3$) points in 
the direction from $L_i$ to $L_j$.
Let $\sgn(F)$ be the product of the signs of all clasp arcs in $F$.

Let $n_i$ be the positive unit normal vector field to $F_i$.
(The co-orientation of $F_i$ is determined by the orientations of $F_i$ and $S^3$.)
We may extend $n_i$ to a smooth vector field $v_i$ on $S^3$ with support in 
a small open neighborhood $U_i$ of $F_i$ such that $U_i$ admits a deformation 
retraction $r_t$ onto $F_i$ and $v_i(x)$ is a positive scalar multiple of $v_i\big(r_1(x)\big)$ 
for each $x\in U_i$.%
\footnote{In more detail, $F_i$ extends to a larger embedded smooth surface $F_i^+$, where 
$\overline{F_i^+\but F_i}$ is the image of a smooth embedding $e_i\:L_i\x [0,1]\to S^3$ 
such that $e_i(x,1)=x$ for $x\in L_i$.
A closed neighborhood $\bar U_i$ of $F_i$ can then be parameterized as the image of 
a smooth embedding $h_i\:F_i^+\x [-1,1]\to S^3$ such that $h_i(x,0)=x$ for $x\in F_i^+$.
Now $v_i$ is defined to be zero outside $\bar U_i$ and by 
$v_i\big(h_i(x,t)\big)=(1-|t|)\pi_i(x)\bar n_i(x)$, where $\bar n_i$ is the positive unit normal vector 
field to $F_i^+$ and $\pi_i\:F_i^+\to[0,1]$ is defined by $\pi_i\big(e_i(y,s)\big)=s$ and 
$\pi_i(x)=1$ if $x\in F$.}

Given a loop $\alpha$ in $|F|$ and signs $\epsilon_1,\dots,\epsilon_n\in\{-1,1\}$, 
let $\alpha^{\epsilon_1\dots\epsilon_n}$ be the pushoff of $\alpha$ (away from $|F|$) along the vector
field $\epsilon_1v_1+\dots+\epsilon_nv_n$.
It is easy to see that the assignment $[\alpha]\mapsto[\alpha^{\epsilon_1\dots\epsilon_n}]$ yields
a well-defined map $H_1(|F|)\to H_1(S^3\but |F|)$, whose value on the homology class of a loop
which is disjoint from $F_i$ does not depend on $\epsilon_i$.
In particular, $\lk(\alpha^{\epsilon_1\dots\epsilon_n},\beta)$ is defined for any $\alpha,\beta\in H_1(|F|)$.

\begin{example} \label{self-linking}
Let $\alpha$ be a smoothly embedded loop in $|F|$ which passes along some of the clasp arcs
$J_{ijk}$, where $i<j$, in the direction of their orientation and is disjoint from the other ones.
Let $\delta_i=(\overbrace{+\dots+}^i\overbrace{-\dots-}^{n-i})$.
Then $\lk(\alpha^{\delta_i},\,\alpha)$ is the same for all $i=0,\dots,n$.

Indeed, given an $i>0$, let us show that $\lk(\alpha^{\delta_i},\alpha)=\lk(\alpha^{\delta_{i-1}},\alpha)$.
Let $\nu_i$ be a normal vector field to $\alpha$ which restricts to the positive unit normal vector field 
to $\alpha\cap F_i$ within $F_i$ and vanishes outside a small neighborhood of $\alpha\cap F_i$.
Let $\alpha_{i,t}$, where $t\in[-1,1]$, be the pushoff of $\alpha$ along the vector field 
$tv_i+(1-|t|)\nu_i+\sum_{j<i}v_j-\sum_{j>i}v_j$.
Thus $\alpha_{i,1}=\alpha^{\delta_i}$ and $\alpha_{i,-1}=\alpha^{\delta_{i-1}}$.
Clearly $\nu_i|_{J_{jik}}$ equals $\pm v_j|_{J_{jik}}$ for each clasp arc $J_{jik}$, 
the sign being determined by the orientations.
Whenever $\alpha$ passes along a clasp arc $J_{jik}$, by our hypothesis their orientations agree 
if $j<i$ and disagree if $j>i$.
Then it is not hard to check that $\nu_i|_{J_{jik}}$ equals $v_j|_{J_{jik}}$ if $j<i$ and $-v_j|_{J_{jik}}$ if $j>i$.
It follows that $\alpha_{i,t}$ is disjoint from $\alpha$ for all $t\in [-1,1]$.
\end{example}

\begin{theorem} \label{cimasoni} 
\cite{Coo2}, \cite{Ci0}, \cite{DMO} Let $F=(F_1,\dots, F_n)$ be a C-complex for a colored smooth link 
$L=(L_1,\dots,L_n)$.
Let $\alpha_1,\dots,\alpha_r$ be a basis for $H_1(|F|)$, and let $A^{\epsilon_1,\dots,\epsilon_n}$
be the Seifert matrix $\big(\lk(\alpha_i^{\epsilon_1,\dots,\epsilon_n},\,\alpha_j)\big)$ for this basis.
Then 
\[\Omega_L(x_1,\dots,x_n)=\sgn(F)\,\prod_{i=1}^n(x_i-x_i^{-1})^{\chi(\bigcup_{j\ne i}F_j)-1}
\det\Big(-\sum_{\epsilon_1,\dots,\epsilon_n=\pm1}\epsilon_1\cdots\epsilon_n\ x_1^{\epsilon_1}\cdots x_n^{\epsilon_n}\,
A^{\epsilon_1\dots\epsilon_n}\Big).\]
\end{theorem}

\begin{corollary} \label{4.1'}
Let $L$ be a colored link with $m$ components.
 
(a) $\Omega_L(x_1,\dots,x_n)=\Omega_\Lambda(-x_1^{-1},\dots,-x_n^{-1})$.

(b) For $m>1$ the total degree of every nonzero term of $\Omega_L$ is congruent $\bmod 2$ to $m$.

(c) For $m>1$ the exponent of $x_i$ in every nonzero term of $\Omega_L$ is congruent $\bmod 2$ 
to $l_i+m_i$, where $m_i$ is the number of components of color $i$ and $l_i=\lk(L_i,\,L\but L_i)$, 
where $L_i$ is the sublink of color $i$.
\end{corollary}

See Lemma \ref{4.1} for an alternative proof.

\begin{proof}[Proof. (a)]
We have $\lk(\alpha_i^{-\epsilon_1,\dots,-\epsilon_n},\,\alpha_j)=\lk(\alpha_i,\,\alpha_j^{\epsilon_1,\dots,\epsilon_n})
=\lk(\alpha_j^{\epsilon_1,\dots,\epsilon_n},\,\alpha_i)$, so
$A^{-\epsilon_1,\dots,-\epsilon_n}$ is the transpose of $A^{\epsilon_1\dots\epsilon_n}$.
Hence $\epsilon_1\cdots\epsilon_n\,(-x_1^{-1})^{\epsilon_1}\cdots (-x_n^{-1})^{\epsilon_n}\,
A^{\epsilon_1,\dots,\epsilon_n}$
is the transpose of $(-\epsilon_1)\cdots(-\epsilon_n)\,x_1^{-\epsilon_1}\cdots x_n^{-\epsilon_n}\,
A^{-\epsilon_1,\dots,-\epsilon_n}$.
\end{proof}

\begin{proof}[(c)]
Let $F=(F_1,\dots,F_n)$ be a C-complex for the colored link $L=(L_1,\dots,L_n)$.
By Theorem \ref{cimasoni} the exponent $e_k$ of $x_k$ in any given nonzero term of $\Omega_L(x_1,\dots,x_n)$ 
has the same parity as $\rk H_1(|F|)+\chi(F^k)-1$, where $F^k=\bigcup_{i\ne k}F_i$.
It is easy to see, arguing by induction and using the Mayer--Vietoris sequence, that $H_2(|F|)=0$.
Since $|F|$ is connected, we get that $\rk H_1(|F|)=1-\chi(|F|)$.
Hence $e_k$ has the same parity as $\chi(|F|)-\chi(F^k)$.
Since $F_i\cap F_i\cap F_k=\emptyset$ for pairwise distinct $i$, $j$ and $k$, we have
$\chi(|F|)=\sum_i\chi(F_i)-\sum_{i<j}\chi(F_i\cap F_j)$.
The parity of $\chi(F_i)$ is the same as that of the number $m_i$ of components in $\partial F_i$.
If $F_i\cap F_j$ consists of $p$ positive and $n$ negative clasps, then $p+n=\chi(F_i\cap F_j)$
and $p-n=\lk(L_i,L_j)$, so the parity of $\chi(F_i\cap F_j)$ is the same as that of $\lk(L_i,L_j)$.
Thus we get $\chi(|F|)\equiv\sum_i m_i+\sum_{i<j}\lk(L_i,L_j)\pmod2$.
Similarly $\chi(F^k)\equiv\sum_{i\ne k}m_i+\sum_{i<j,\,i\ne k\ne j}\lk(L_i,L_j)\pmod2$.
Therefore $e_k$ has the same parity as $m_k+\sum_{i\ne k}\lk(L_k,L_i)=m_k+\lk(L_k,\,L\but L_k)$.
\end{proof}

\begin{proof}[(b)] This follows immediately from (c).
\end{proof}

\begin{example} \label{kauffman} Let $L$ be a smooth link and $F$ its connected Seifert surface.
Then $A^\epsilon=\big(\lk(\alpha_i^\epsilon,\,\alpha_j)\big)$.
We have $\lk(\alpha_i^-,\,\alpha_j)=\lk(\alpha_i,\alpha_j^+)=\lk(\alpha_j^+,\alpha_i)$, so $A^-$ is
the transpose of $A^+$.
Thus $\Omega_L(x)=(x-x^{-1})^{-1}\det(-xA^++x^{-1}A^-)$.
Therefore $\nabla_L(x-x^{-1})=(x-x^{-1})\Omega_L(x)=\det(-xV+x^{-1}V^T)$, where $V=A^+=\big(\lk(\alpha_i^+,\,\alpha_j)\big)$ 
is the usual Seifert matrix of $F$.
This formula for the Conway polynomial $\nabla_L$ appears in Cimasoni's paper \cite{Ci} and is originally due to 
Kauffman \cite{Kau} (see also \cite{Kau2}).%
\footnote{But Kauffman writes it differently, because his ``linking number'' is our $-\lk$ 
(compare \cite{Kau2}) and consequently his ``Seifert matrix'' is our $-V$.}

Let us note that since $\nabla_K(z)$ involves only even powers of $z$, we have $\nabla_K(x-x^{-1})=\nabla_K(x^{-1}-x)$,
and so one can also write $\nabla_L(x-x^{-1})=\det(-x^{-1}V+xV^T)$.
On the other hand, the size of $V$ is $2g+(m-1)$, where $m$ is the number of components of $L$ and $g$ is the genus of 
the closed surface obtained by attaching disks to $F$.
Thus we may also write $\nabla_L(x-x^{-1})=(-1)^{m-1}\det(xV-x^{-1}V^T)$.

For a knot $K$, the Alexander polynomial $\Delta_K(t)$, $t=x^2$, which is well-defined up to units of $\Z[t]$,
is represented by $\nabla_K(x-x^{-1})=\det(x^{-1}V-xV^T)=x^{-2g}\det(V-x^2V^T)$, which is the same as $\det(V-tV^T)$ 
up to units of $\Z[t]$.
This formula is originally due to Seifert \cite{Se}.
\end{example}

\begin{example} \label{kauffman-sign} Let us discuss sign conventions involved in the formula
$\nabla_L(x-x^{-1})=\det(-xV+x^{-1}V^T)$.
Let $K$ be the unknot and let $K'$ be its pushoff such that $\lk(K,K')=+1$.
Let $F$ be the obvious annulus cobounded by $K$ and $-K'$.
Let $Q$ be a pushoff of $K$ within $F$.
As usual, $Q^+$ is a pushoff of $Q$ off $F$ in the positive direction.
Then $V$ is the $1\x 1$ matrix whose only entry is $\lk(Q^+,Q)=\lk(K,Q)=1$.
Hence $\nabla_L(z)=-z$, where $L=(K,-K')$, which is consistent with $\lk(L)=-1$ being the coefficient
of the linear term of $\nabla_L(z)$.
Let us note that the co-orientation of $F$ is determined by various sign conventions; but if it is reversed, 
$Q^+$ becomes $Q^-$, and we have $\lk(Q^-,Q)=\lk(Q,Q^+)=\lk(Q^+,Q)$ --- so this does not affect the sign of $\nabla_L$.

One can also compute $\nabla_L$ from the skein relation 
$\nabla_{\includegraphics[width=0.4cm]{1+p.pdf}}(z)-\nabla_{\includegraphics[width=0.4cm]{1-p.pdf}}(z)
=z\,\nabla_{\includegraphics[width=0.4cm]{1op.pdf}}$.
If we agree that $\lk(\raisebox{-0.15cm}{$\includegraphics[width=0.9cm]{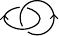}$})=+1$, 
then \raisebox{-0.15cm}{$\includegraphics[width=0.9cm]{Hopfp.pdf}$}
is a picture of $(K,K')$; and it follows that the skein relation can be applied so that
\raisebox{-0.1cm}{$\includegraphics[width=0.5cm]{1-p.pdf}$} represents $L$, 
\raisebox{-0.1cm}{$\includegraphics[width=0.5cm]{1+p.pdf}$} represents the $2$-component unlink,
and \raisebox{-0.1cm}{$\includegraphics[width=0.5cm]{1op.pdf}$} represents the unknot.
This again yields $\nabla_L(z)=-z$.

To summarize, as long as the linking number is defined in a way consistent with the skein relation 
$\nabla_{\includegraphics[width=0.4cm]{1+p.pdf}}(z)-\nabla_{\includegraphics[width=0.4cm]{1-p.pdf}}(z)
=z\,\nabla_{\includegraphics[width=0.4cm]{1op.pdf}}$, that is, so that 
$\lk(\raisebox{-0.1cm}{$\includegraphics[width=0.5cm]{1+p.pdf}$})-\lk(\raisebox{-0.1cm}{$\includegraphics[width=0.5cm]{1-p.pdf}$})=1$
for crossing changes involving both components, or equivalently $\dfrac{d\nabla_L(z)}{dz}\Big|_{z=0}=\lk(L)$, the sign of 
$\det(-xV+x^{-1}V^T)$ does not depend on any further sign conventions.

See also \cite{Tr2}*{\S2} concerning sign conventions for $\Omega_L$.
\end{example}

\begin{example} \label{cooper}
Let $L=(K_1,K_2)$ be a $2$-component smooth link and $F=(F_1,F_2)$ be a C-complex for $L$.
Suppose that $F$ contains $c$ clasps, of which $p$ are positive and $n$ are negative.
Thus $\sgn(F)=(-1)^n$ and $p+n=c$.
Also $p-n=\lk(L)$, whence $n=\big(c-\lk(L)\big)/2$.
On the other hand, from the Mayer--Vietoris exact sequence 
\[0\to H_1(F_1)\oplus H_1(F_2)\to H_1(F_1\cup F_2)\to \tilde H_0(F_1\cap F_2)\to 0\]
the size $r$ of the Seifert matrix equals $\rk H_1(F_1)+\rk H_1(F_2)+c-1$,
which is of the same parity as $\lk(L)+1$.
Hence $\det(-M)=(-1)^{\lk(L)+1}\det M$ for a square matrix of size $r$.
Therefore
\[\Omega_L(x,y)=(-1)^{(c+\lk(L)+2)/2}\,(x-x^{-1})^{\chi(F_2)-1}(y-y^{-1})^{\chi(F_1)-1}
\det\Big(\sum_{\epsilon,\delta=\pm1}\epsilon\delta\ x^\epsilon y^\delta\,A^{\epsilon\delta}\Big).\]
This is how the sign of $\Omega_L$ is described in \cite{Coo1}*{Proposition 7.4}.
The description of Theorem \ref{cimasoni}(b) appears in \cite{Ci0}.
\end{example}

\section{Factorization: linking number 1}

Now let $L=(K_1,K_2)$ be a $2$-component smooth link with $\lk(L)=1$.
Then by Lemma \ref{C-complex} $L$ has a C-complex $F=(F_1,F_2)$ with just one clasp.
Let $\alpha_1,\dots,\alpha_{2g_1}$ be a symplectic basis of $H_1(F_1)$ and
$\alpha_{2g_1+1},\dots,\alpha_{2g_1+2g_2}$ be a symplectic basis of $H_1(F_2)$.
Let us note that together they form a basis of $H_1(|F|)$.
The Seifert matrices $A^{\epsilon\delta}=\big(\lk(\alpha_i^{\epsilon\delta},\,\alpha_j)\big)$
are of the following block form:
\[A^{++}=\begin{pmatrix} V_1 & L\\ L^T & V_2 \end{pmatrix},\ \
A^{+-}=\begin{pmatrix} V_1 & L\\ L^T & V_2^T \end{pmatrix},\ \
A^{-+}=\begin{pmatrix} V_1^T & L\\ L^T & V_2 \end{pmatrix},\ \
A^{--}=\begin{pmatrix} V_1^T & L\\ L^T & V_2^T \end{pmatrix},\]
where $V_i$ is the Seifert matrix for $F_i$ with respect to the chosen basis in $H_1(F_i)$, and $L$ is the matrix 
$\big(\lk(\alpha_i,\alpha_j)\big)_{i\le 2g_1<j}$.
Let $M_1=xV_1-x^{-1}V_1^T$ and $M_2=yV_2-y^{-1}V_2^T$.

\begin{theorem} \label{factorization0} In the above notation,
\[\Omega_L(x,y)=\nabla_{K_1}(x-x^{-1})\,\nabla_{K_2}(y-y^{-1})\,\det\Big(I-(x-x^{-1})(y-y^{-1})M_1^{-1}LM_2^{-1}L^T\Big).\]
\end{theorem}

\begin{proof} By Example \ref{cooper}
\begin{multline*}\Omega_L(x,y)=(x-x^{-1})^{-2g_2}\,(y-y^{-1})^{-2g_1}\big|xyA^{++}-xy^{-1}A^{+-}-x^{-1}yA^{-+}+x^{-1}y^{-1}A^{--}\big|\\
=(x-x^{-1})^{-2g_2}\,(y-y^{-1})^{-2g_1}\left|\begin{matrix}
(y-y^{-1})(xV_1-x^{-1}V_1^T) & (x-x^{-1})(y-y^{-1})L^T \\[3pt] 
(x-x^{-1})(y-y^{-1})L & (x-x^{-1})(yV_2-y^{-1}V_2^T)\end{matrix}\right|\\
=\left|\begin{matrix}
M_1 & (x-x^{-1})L^T \\[3pt]
(y-y^{-1})L & M_2
\end{matrix}\right|.
\end{multline*}
By Lemma \ref{schur} 
\begin{multline*}
\Omega_L(x,y)=\det(M_2)\,\det\big(M_1-(x-x^{-1})(y-y^{-1})LM_2^{-1}L^T\big)\\
=\det(M_1)\,\det(M_2)\,\det\big(I-(x-x^{-1})(y-y^{-1})M_1^{-1}LM_2^{-1}L^T\big).
\end{multline*}
Now the desired formula follows from Example \ref{kauffman}.
\end{proof}

\begin{lemma} \label{schur} For a block matrix $M=\begin{pmatrix} A & C\\ C' & B \end{pmatrix}$
where $A$ and $B$ are invertible%
\footnote{Let us recall that a matrix $A$ over a ring is invertible if and only if $\det A$
is invertible (see Lemma \ref{adjugate}).
This assertion is more widely known for matrices over a field, which case will in fact suffice for
the purposes of our applications of Lemma \ref{schur}.}
square matrices $\det(M)=\det(A)\,\det(B-C'A^{-1}C)=\det(B)\,\det(A-CB^{-1}C')$.
\end{lemma}

This is well-known. $A-CB^{-1}C'$ and $B-C'A^{-1}C$ are called the {\it Schur complements}.

\begin{proof} The first assertion follows from
$\begin{pmatrix} A & C\\ C' & B \end{pmatrix}
\begin{pmatrix} I & -A^{-1}C\\ 0 & I \end{pmatrix}
=\begin{pmatrix} A & 0\\ C' & B-C'A^{-1}C \end{pmatrix}$
and the second from
$\begin{pmatrix} I & -CB^{-1}\\ 0 & I \end{pmatrix}
\begin{pmatrix} A & C\\ C' & B \end{pmatrix}
=\begin{pmatrix} A-CB^{-1}C' & 0\\ C' & B \end{pmatrix}$.
\end{proof}

\begin{remark} \label{factorization0'} From the proof of Theorem \ref{factorization0} we also have
\begin{align*}\Omega_L(x,y)&=\nabla_{K_2}(y-y^{-1})\,\det\big(M_1-(x-x^{-1})(y-y^{-1})LM_2^{-1}L^T\big)\\
&=\nabla_{K_1}(x-x^{-1})\,\det\big(M_2-(x-x^{-1})(y-y^{-1})L^TM_1^{-1}L\big).
\end{align*}
\end{remark}

Let $A_1,\dots,A_{2g_1}\subset F_1$ be a symplectic basis and $B_1,\dots,B_{2g_1}\subset S^3\but F_1$ 
the associated collection of meridians (see \S\ref{cochran-pairing-section}).

\begin{corollary} \label{factorization1}
\[\Omega_L(x,y)=\nabla_{K_1}(x-x^{-1})\,\nabla_{K_2}(y-y^{-1})\,\det(I-SR),\] 
where $R=\Big(\langle A_i^+,A_j\rangle_{F_2}-\lk(A_i^+,A_j)\Big)_{1\le i,j\le 2g_1}$ 
and $S=\Big(\langle B_i,B_{j+}\rangle_{F_1}\Big)_{1\le i,j\le 2g_1}$.
\end{corollary}

\begin{proof} By Proposition \ref{PY1}(b) $S=(x-x^{-1})(-M_1)^{-1}$.
By Theorem \ref{PY2}(b) (with $y$ in place of $x$) $R=(y-y^{-1})L(-M_2)^{-1}L^T$.
Now the assertion follows from Theorem \ref{factorization0}.
\end{proof}

It is convenient to use the notation $C_L(s)=\sum_{k=0}^\infty\alpha_{1,\,2k-1} s^k$.

\begin{theorem} \label{factorization2}
For a smooth link $L=(K_1,K_2)$ with $\lk(L)=1$ we have
\[C_L(s)=(-1)^{g_1}\sum_{n=1}^\infty(-1)^n\sum_{i=1}^{g_1}\beta_{F_2}^{n-1,\,n}(A_{2i},A_{2i-1})\,s^{n-1},\]
where $F_1$, $F_2$ are Seifert surfaces for $K_1$ and $K_2$ which intersect along a single clasp arc, 
and $A_1,\dots,A_{2g_1}$ is a symplectic basis for $F_1$, disjoint from the clasp.
\end{theorem}

\begin{proof}
Since multiplication by $\nabla_K(x_i-x_i^{-1})$ commutes with the endomorphism $\bar\Theta_L$ (see \S\ref{reduced})
$\bar\varnabla_L(u,v)$ can be uniquely presented in the form 
$\bar\varnabla(u,v)=\sum_{i=0}^\infty u^{2i}P_{2i}(v^2)+w\sum_{i=0}^\infty u^{2i+1}vP_{2i+1}(v^2)$ and clearly $C_{\bar L}(v^2)=P_1(v^2)$.
Similarly $\varnabla_L(u,v)/\nabla_{K_2}(v)$ can be uniquely presented in the form
$\varnabla_L(u,v)/\nabla_{K_2}(v)=\sum_{i=0}^\infty u^{2i}Q_{2i}(v^2)+w\sum_{i=0}^\infty u^{2i+1}vQ_{2i+1}(v^2)$.
Since $\nabla_{K_1}(u)$ is of the form $1+c_1u^2+c_2u^4+\dots$, it is easy to see that $P_1(v^2)=Q_1(v^2)$.

From Remark \ref{factorization0'} we get $\Omega_L(x,y)=\nabla_{K_2}(v)\,\det(M_1+uR)$, where
$R=vL(-M_2)^{-1}L^T$.
By Theorem \ref{PY2}(b) $R=\big(\left<A_i^+,A_j\right>_{F_2}-\lk(A_i^+,A_j)\big)$.
Therefore by Theorem \ref{cochran expansion}(b) $R=(r_{ij})$, where
\[r_{ij}=\sum_{n=1}^\infty(-1)^n\Big(y\beta_{F_2}^{n,n}(A_j,A_i^+)-y^{-1}\beta_{F_2}^{n,n}(A_i^+,A_j)\Big)v^{2n-1}.\]
Here we may replace $A_i^+$ by $A_i$ since $\beta_{F_2}^{n,n}(A,B)$ depends only on $[A],[B]\in H_1(S^3\but F_2)$
for $n\ge 1$.

Let us recall that $M_1=xV_1-x^{-1}V_1^T$, where $V_1=\big(\lk(\alpha_i^+,\alpha_j)\big)$ is the Seifert matrix of $F_1$ with respect 
to the basis $\alpha_1,\dots,\alpha_{2g_1}$.
By Corollary \ref{seifert-intersection2} $V_1-V_1^T$ is the intersection matrix $J$ of $H_1(F_1)$ with respect to the chosen basis, 
that is, the block diagonal matrix where all diagonal blocks are $\begin{pmatrix} 0 & 1\\ -1 & 0\end{pmatrix}$.
Let $S$ be the symmetric matrix whose entries on and above the diagonal equal those of $V_1$, and let $U=S-V_1$.
Then $J=V_1-V_1^T=(S-U)-(S-U^T)=U^T-U$.
Now $U$ has non-zero entries only below the diagonal, and $U^T$ only above the diagonal.
Hence they are the block diagonal matrices where all diagonal blocks are $\begin{pmatrix} 0 & 1\\ 0 & 0\end{pmatrix}$
and $\begin{pmatrix} 0 & 0\\ 1 & 0\end{pmatrix}$, respectively.
Now $M_1=xV_1-x^{-1}V_1^T=x(S-U)-x^{-1}(S-U^T)=(x-x^{-1})S-xU+x^{-1}U^T=uS+H$,
where $H$ is the block diagonal matrix where all diagonal blocks are $\begin{pmatrix} 0 & x^{-1}\\ -x & 0\end{pmatrix}$.

To summarize, $\Omega_L(x,y)=\nabla_{K_2}(v)\,\det(uS+H+uR)$.
Thus the power series $C_L(v^2)=P_1(v^2)=Q_1(v^2)$ can be computed by evaluating $\det(uS+H+uR)$, 
putting it in the standard form and collecting all terms that are divisible by $uv$ but not divisible by $u^2$.
We have $\det(a_{ij})=\sum_{\sigma\in S_{2g_1}}\sgn(\sigma)a_{1,\sigma(1)}\cdots a_{2g_1,\sigma(2g_1)}$.
For a summand in this sum to be not divisible by $u^2$, at most one of its factors can be divisible by $u$.
Thus we can take at most one factor from $uS$ or $uR$, and all other factors must be taken from $H$.
This condition is satisfied only for one permutation $\sigma$, namely, the one which gives a non-zero contribution
to $\det H$.
On the other hand, if we take all the factors from $H$ and $uS$, then the resulting summand will involve no occurrences 
of $y$ (since $uS=M_1-H$ and neither $M_1$ nor $H$ involves any occurrences of $y$) and so will give no contribution 
to $Q_1(v^2)$.
Thus we must take precisely one factor from $uR$ and all the remaining factors from $H$.
There are $2g_1$ ways of doing so, and here is the sum of the resulting summands:
$(-1)^{g_1-1}\sum_{i=1}^{g_1}\big(-xur_{2i-1,\,2i}+x^{-1}ur_{2i,2i-1}\big)$.
This can be rewritten as
\[(-1)^{g_1}u\sum_{i=1}^{g_1}\sum_{n=1}^\infty(-1)^n
\Big((xy+x^{-1}y^{-1})\beta_{F_2}^{n,n}(A_{2i},A_{2i-1})-(xy^{-1}+x^{-1}y)\beta_{F_2}^{n,n}(A_{2i-1},A_{2i})\Big)v^{2n-1}.\]
Now we apply the identity $xy+x^{-1}y^{-1}=uv+x^{-1}y+xy^{-1}$, and upon discarding summands that are divisible by $u^2$, we get
\[(-1)^{g_1}(xy^{-1}+x^{-1}y)u\sum_{i=1}^{g_1}\sum_{n=1}^\infty(-1)^n
\Big(\beta_{F_2}^{n,n}(A_{2i},A_{2i-1})-\beta_{F_2}^{n,n}(A_{2i-1},A_{2i})\Big)v^{2n-1}.\]
This can now be equated with $(xy^{-1}+x^{-1}y)uvC_L(v^2)$, which yields the desired formula,
taking into account that 
$\beta_{F_2}^{n,n}(A,B)-\beta_{F_2}^{n,n}(B,A)=\beta_{F_2}^{n-1,\,n}(A,B)$ by Lemma \ref{pascal}(d).
\end{proof}

\section{Factorization: linking number 0}

The following is well-known; (b) is called the Sherman--Morrison--Woodbury formula.

\begin{lemma} \label{woodbury} For a block matrix $M=\begin{pmatrix} A & C\\ C' & B \end{pmatrix}$
where $A$ and $B$ are invertible square matrices let us consider the Schur complements
$\hat A=A-CB^{-1}C'$ and $\hat B=B-C'A^{-1}C$ and assume that they are invertible
(or equivalently $M$ is invertible, see Lemma \ref{schur}).
Then

(a) $\hat A^{-1}CB^{-1}=A^{-1}C\hat B^{-1}$ and $\hat B^{-1}C'A^{-1}=B^{-1}C'\hat A^{-1}$.
\smallskip

(b) $\hat A^{-1}-A^{-1}=A^{-1}C\hat B^{-1}C'A^{-1}$ and $\hat B^{-1}-B^{-1}=B^{-1}C'\hat A^{-1}CB^{-1}$.
\end{lemma}

The second assertions of (a) and (b) reduce to the first ones by renaming the matrices.

\begin{proof}[Proof. (a)] We have
\[CB^{-1}\hat B=CB^{-1}(B-C'A^{-1}C)=C-CB^{-1}C'A^{-1}C=(A-CB^{-1}C')A^{-1}C=\hat AA^{-1}C.\]
Now multiply both sides by $\hat A^{-1}$ on the left and by $\hat B^{-1}$ on the right.
\end{proof}

\begin{proof}[(b)] By (a)
$A^{-1}C\hat B^{-1}C'A^{-1}=\hat A^{-1}CB^{-1}C'A^{-1}=\hat A^{-1}(A-\hat A)A^{-1}=\hat A^{-1}-A^{-1}.$
\end{proof}

\begin{lemma} \label{3x3} Let us consider a block matrix
\[N=\begin{pmatrix}
A & C & \Gamma \\
C' & B & \Delta \\
\Gamma' & \Delta' & Q
\end{pmatrix},\]
where $A$, $B$ and $Q$ are square matrices, with $A$ and $B$ invertible.
Assume further that $M:=\begin{pmatrix}
A & C \\
C' & B \\
\end{pmatrix}$
is invertible, or equivalently (see Lemma \ref{schur})
the Schur complements $\hat A:=A-CB^{-1}C'$ and $\hat B:=B-C'A^{-1}C$ are invertible.
Then \[\det(N)=\det(M)\,\det\Big(Q-\Gamma'\hat A^{-1}\Gamma-\Delta'\hat B^{-1}\Delta
+\Gamma'\hat A^{-1}CB^{-1}\Delta+\Delta'\hat B^{-1}C'A^{-1}\Gamma\Big).\] 
\end{lemma}

\begin{proof} Let us find matrices $X$ and $Y$ of the same sizes as $\Gamma$ and $\Delta$ such that
\[\begin{pmatrix}
A & C & \Gamma \\
C' & B & \Delta \\
\Gamma' & \Delta' & Q
\end{pmatrix}
\begin{pmatrix}
I & 0 & X \\
0 & I & Y \\
0 & 0 & I
\end{pmatrix}
=
\begin{pmatrix}
A & C & 0 \\
C' & B & 0 \\
Z_1 & Z_2 & Q'
\end{pmatrix}\]
for some matrices $Z_1$, $Z_2$ and $Q'$.
Thus we have the system of matrix equations
\[\begin{cases}
AX+CY=-\Gamma\\
C'X+BY=-\Delta.
\end{cases}\]
Multiplying the second equation by $CB^{-1}$ on the left, we get $CB^{-1}C'X+CY=-CB^{-1}\Delta$.
Subtracting this from the first equation, we get $(A-CB^{-1}C')X=CB^{-1}\Delta-\Gamma$,
that is, $\hat AX=CB^{-1}\Delta-\Gamma$.%
\footnote{Two remarks due to M. Il'insky: (i) $CB^{-1}\Delta-\Gamma$ is a Schur complement
for the upper right $2\x 2$ submatrix of $N$; (ii) by symmetry, from 
$\hat AX=CB^{-1}\Delta-\Gamma$ one can infer that $\hat BY=C'A^{-1}\Gamma-\Delta$.}
Since $M$ and $A$ are invertible, so is $\hat A$, so we get
$X=\hat A^{-1}CB^{-1}\Delta-\hat A^{-1}\Gamma$.
Substituting for $X$ in the second equation, we find that
$Y=-B^{-1}C'X-B^{-1}\Delta=B^{-1}C'\hat A^{-1}\Gamma-B^{-1}C'\hat A^{-1}CB^{-1}\Delta-B^{-1}\Delta$.
Using Lemma \ref{woodbury}(a,b) this can be rewritten as 
$Y=\hat B^{-1}C'A^{-1}\Gamma-\hat B^{-1}\Delta$.
Now $Q'=\Gamma'X+\Delta'Y+Q$, and by substituting for $X$ and $Y$ we get the desired expression.
\end{proof}

\begin{lemma}\cite{AADG} \label{C-complex2}
Let $L=(K_1,K_2)$ be a PL link with $\lk(L)=0$.
Then there exist Seifert surfaces $F_i$ for $K_i$, intersecting transversely and such that 
$F_1\cap F_2$ consists of precisely two clasp arcs.
\end{lemma}

For variety we include a proof based on the arguments of Cooper \cite{Coo1}*{Lemma 3.1}.

\begin{proof} Since $\lk(L)=0$, each $K_i$ bounds an embedded connected surface $F_i$ in $S^3\but K_{3-i}$.
Without loss of generality $F_1$ and $F_2$ intersect transversely and their intersection is nonempty.
If $C:=F_1\cap F_2$ contains more than one connected component, then some two of them are adjacent on
$F_1$ in the sense that there exists an arc $J$ in $F_1$ whose two endpoints lie in these two 
components and which is otherwise disjoint from $F_2$ and from $K_1$.
By attaching to $F_2$ a tube running along $J$ we will make these two components into one.
By continuing this process we may assume that $C$ is connected.

Next let us join $C$ to $K_1$ by an arc $J_1$ in $F_1$ which is otherwise disjoint from $F_2$ 
and from $K_1$.
By pushing a finger out of $F_2$ along $J_1$ we will make $C$ into an arc with both endpoints in $K_1$.
Finally, let us join an interior point of $C$ to $K_2$ by an arc $J_2$ in $F_2$ which is otherwise 
disjoint from $F_1$ and from $K_2$.
By pushing a finger out of $F_1$ along $J_2$ we will make $C$ into a pair of arcs each with 
one endpoint in $K_1$ and another in $K_2$.
\end{proof}

Let $L=(K_1,K_2)$ be a $2$-component smooth link with $\lk(L)=0$, and 
let $F=(F_1,F_2)$ be a C-complex for $L$ with precisely two clasp arcs $J_+$ and $J_-$.
Let $J_i\subset K_i$ be the arc joining the two points $K_i\cap F_{3-i}$ and lying on the positive side of $F_{3-i}$.
Then $\Pi_i:=J_i\cup J_+\cup J_-$ is an arc lying in $F_i$.
A regular neighborhood of $\Pi_i$ in $F_i$ is a disk $D_i$ which meets $K_i$ in an arc.
Let $D=D_1\cup D_2$ and let $S_i=\overline{F_i\but D_i}$.
Each $S_i$ meets $D$ in an arc and $S_1\cap S_2=\emptyset$, so 
$H_1(|F|)\simeq H_1(S_1)\oplus H_1(S_2\cup D)\simeq H_1(S_1)\oplus H_1(S_2)\oplus H_1(D)$.
Clearly each $S_i$ is a deformation retract of $F_i$, whereas $D$ deformation retracts onto
the knot $J_1\cup J_2\cup J_+\cup J_-$.

If we orient $J_+$ and $J_-$ according to the orientation of $F_1\cap F_2$ (rather than $F_2\cap F_1$), 
then one, say $J_+$, is oriented from $K_1$ to $K_2$ and the other one (that is, $J_-$) in the opposite direction.
Let $Q\subset D$ be a knot which first goes along $J_+$, then along
a pushoff of $\Pi_2$ (relative to $\partial\Pi_2$) into the interior of $D_2$, then along $J_-$, and finally along
a pushoff of $\Pi_1$ (relative to $\partial\Pi_1$) into the interior of $D_1$.
Then the orientation of $Q$ agrees with those of $J_+$ and $J_-$.
Hence by Example \ref{self-linking} $\lk(Q^{--},\,Q)=\lk(Q^{+-},\,Q)=\lk(Q^{++},\,Q)$.
Also $\lk(Q^{-+},\,Q)=\lk(Q,\,Q^{+-})=\lk(Q^{+-},\,Q)$.
Thus $\lk(Q^{\epsilon\delta},\,Q)$ does not depend on $\epsilon$ and $\delta$.
Let $\alpha_1,\dots,\alpha_{2g_1}$ be a symplectic basis of $H_1(S_1)$, let
$\alpha_{2g_1+1},\dots,\alpha_{2g_1+2g_2}$ be a symplectic basis of $H_1(S_2)$, and let
$\alpha_r$, where $r=2g_1+2g_2+1$, be the generator of $H_1(D)$ represented by $Q$.
Thus $\alpha_1,\dots,\alpha_r$ together form a basis of $H_1(|F|)$.
The Seifert matrices $A^{\epsilon\delta}=\big(\lk(\alpha_i^{\epsilon\delta},\,\alpha_j)\big)$
are of the following block form:
\begin{gather*}
A^{++}=\begin{pmatrix} V_1 & L & \Lambda_1^T\\ L^T & V_2 & \Lambda_2^T\\ \Lambda_1 & \Lambda_2 & \beta \end{pmatrix},\ \
A^{+-}=\begin{pmatrix} V_1 & L & \Lambda_1^T\\ L^T & V_2^T & \Lambda_2^T\\ \Lambda_1 & \Lambda_2 & \beta \end{pmatrix},\\[6pt]
A^{-+}=\begin{pmatrix} V_1^T & L & \Lambda_1^T\\ L^T & V_2 & \Lambda_2^T\\ \Lambda_1 & \Lambda_2 & \beta \end{pmatrix},\ \
A^{--}=\begin{pmatrix} V_1^T & L & \Lambda_1^T\\ L^T & V_2^T & \Lambda_2^T\\ \Lambda_1 & \Lambda_2 & \beta \end{pmatrix},
\end{gather*}
where 
\begin{itemize}
\item $V_i$ is the Seifert matrix for $F_i$ with respect to the chosen bases in $H_1(S_i)\simeq H_1(F_i)$;
\item $L$ is the matrix $\big(\lk(\alpha_i,\alpha_j)\big)_{i\le 2g_1<j}$;
\item $\Lambda_1=\big(\lk(\alpha_r,\alpha_1),\dots,\lk(\alpha_r,\alpha_{2g_1})\big)$
and $\Lambda_2=\big(\lk(\alpha_r,\alpha_{2g_1+1}),\dots,\lk(\alpha_r,\alpha_{2g_1+2g_2})\big)$;
\item $\beta=\lk(\alpha_r^{++},\alpha_r)=\lk(Q^{++},Q)$.
\end{itemize}

\begin{lemma} \label{satolevine}
(a) $\beta$ is the Sato--Levine invariant of $L$.

(b) $Q^{++}$ is ambient isotopic in $S^3\but L$ to the oriented intersection of a pair of Seifert surfaces
of $K_1$ and $K_2$.
\end{lemma}

\begin{figure}[h]
\includegraphics[width=9.5cm]{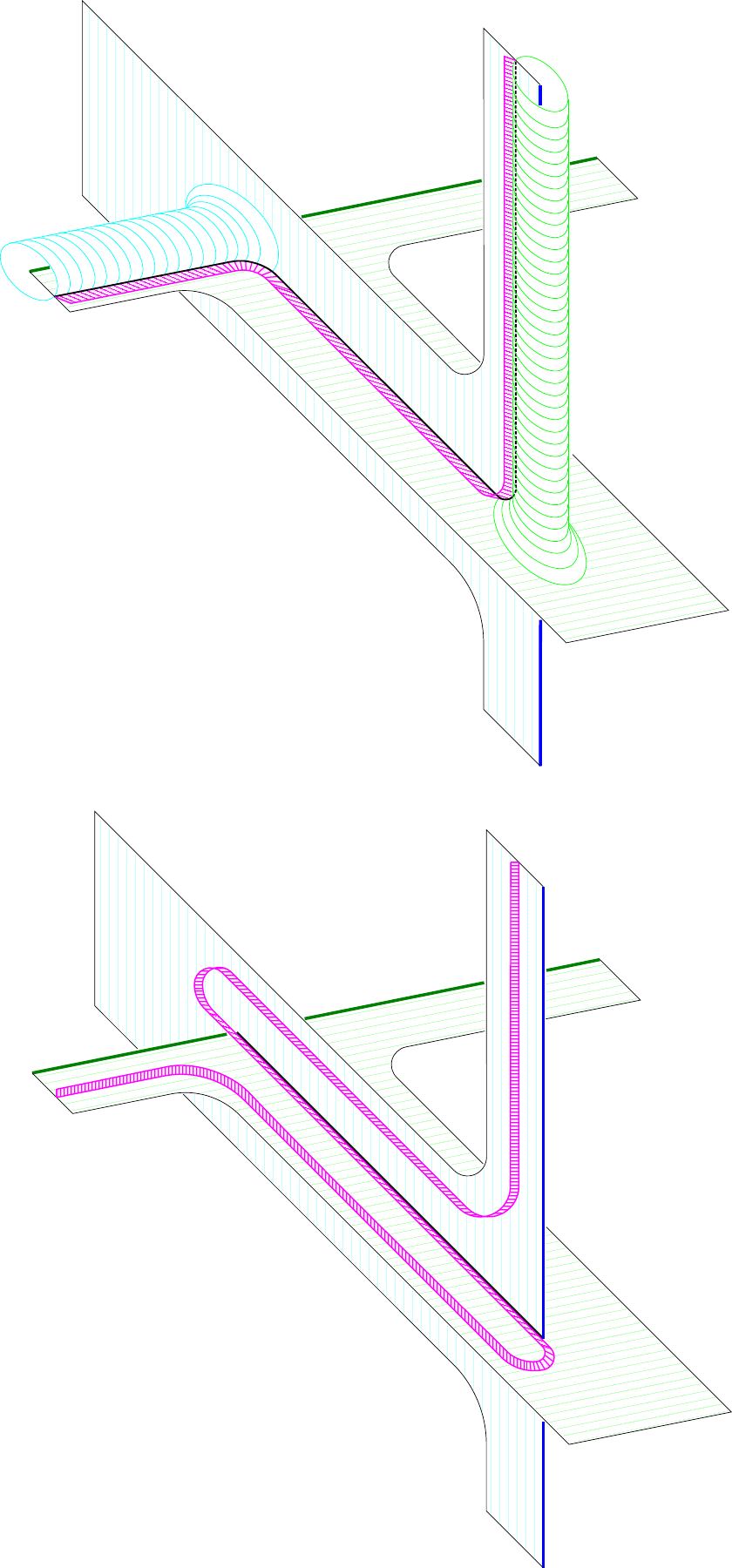}
\caption{Identifying Cochran's derivative.}
\label{cochran-alexander}
\end{figure}

\begin{proof} Let $F_i'$ be obtained by attaching a tube to $F_i$ going along the arc $J_{3-i}$.
Thus $F_i'$ is disjoint from $K_{3-i}$, and the oriented intersection $R:=F_1'\cap F_2'$ is a knot
in $S^3\but L$.
It is easy to see by inspection (see Figure \ref{cochran-alexander}) that
$Q^{++}$ is ambient isotopic in $S^3\but L$ to $R$.
Thus we get (b).
Also it is easy to see by inspection (see again Figure \ref{cochran-alexander}) that
the ribbon cobounded by $Q$ and $Q^{++}$ is ambient isotopic in $S^3$ to the ribbon
cobounded by $R$ and $R^{++}$, where $R^{++}$ is the pushoff of $R$ away from $F_1'$ and $F_2'$
in the direction of the positive co-orientations of both.
Now the Sato--Levine invariant of $L$ by definition equals $\lk(R^{++},R)$, so we get (a).
\end{proof}

We recall the notation $C_L(z)=\sum_{k=0}^\infty\alpha_{1,\,2k-1} z^{2k}$.

\begin{theorem}\label{jin} Let $L=(K_1,K_2)$ be a link with $\lk(L)=0$. 

(a) {\rm (Bailey \cite{Hi}*{7.1})} $\Omega_L(x,y)$ as a Laurent polynomial is divisible by $(x-x^{-1})(y-y^{-1})$.

(b) \cite{Jin}, \cite{M04} $C_L(z)=\sum_{k=0}^\infty(-1)^k\beta^{k+1} z^{2k}$, where $\beta_k$ are Cochran's derived invariants.
\end{theorem}
 
\begin{proof} Let $u=x-x^{-1}$ and $v=y-y^{-1}$.
Let $M_1=xV_1-x^{-1}V_1^T$ and $M_2=yV_2-y^{-1}V_2^T$.
By Example \ref{cooper}
\begin{multline*}\Omega_L(x,y)=u^{-2g_2}v^{-2g_1}\det(xyA^{++}-xy^{-1}A^{+-}-x^{-1}yA^{-+}+x^{-1}y^{-1}A^{--})\\[10pt]
=u^{-2g_2}v^{-2g_1}\left|\begin{matrix}
v(xV_1-x^{-1}V_1^T) & uvL & uv\Lambda_1^T\\[3pt] 
uvL^T & u(yV_2-y^{-1}V_2^T) & uv\Lambda_2^T\\[3pt]
uv\Lambda_1 & uv\Lambda_2 & uv\beta
\end{matrix}\right|
=uv\left|\begin{matrix}
M_1 & uL & u\Lambda_1^T\\[3pt] 
vL^T & M_2 & v\Lambda_2^T\\[3pt]
\Lambda_1 & \Lambda_2 & \beta
\end{matrix}\right|.
\end{multline*}
This completes the proof of (a).

Next, let us consider $M:=\left(\begin{smallmatrix}
M_1 & uL \\[3pt] 
vL^T & M_2 \\[3pt]
\end{smallmatrix}\right)$
and the corresponding Schur complements 
$\hat M_1:=M_1-uvLM_2^{-1}L^T$ and $\hat M_2:=M_2-uvL^TM_1^{-1}L$.
Then by Lemma \ref{3x3} we get
\[\Omega_L(x,y)=uv\beta'\det M,\] where 
$\beta'=\beta-u\Lambda_1\hat M_1^{-1}\Lambda_1^T-v\Lambda_2\hat M_2^{-1}\Lambda_2^T
+uv\Lambda_1\hat M_1^{-1}LM_2^{-1}\Lambda_2^T+uv\Lambda_2\hat M_2^{-1}L^TM_1^{-1}\Lambda_1^T$.
By the proof of Theorem \ref{factorization0}
$\det M=\nabla_{K_1}(u)\,\nabla_{K_2}(v)\,\det\big(I-uvM_1^{-1}LM_2^{-1}L^T\big)$.
Hence \[\bar\Omega_L(x,y)=uv\beta'\det\big(I-uvM_1^{-1}LM_2^{-1}L^T\big).\]
Therefore $\bar\Omega_L(x,y)\equiv uv\beta'\bmod{(u^2)}$.
Then $\bar\Omega_L(x,y)\equiv uv\big(\beta-v\Lambda_2\hat M_2^{-1}\Lambda_2^T\big)\bmod{(u^2)}$.
By Lemma \ref{woodbury}(b) $\hat M_2^{-1}=M_2^{-1}+uvM_2^{-1}L^T\hat M_1^{-1}LM_2^{-1}$.
Hence in fact $\bar\Omega_L(x,y)\equiv uv\big(\beta-v\Lambda_2 M_2^{-1}\Lambda_2^T\big)\bmod{(u^2)}$.
It follows that \[C_L(v)=\beta-v\Lambda_2 M_2^{-1}\Lambda_2^T.\]
By Theorem \ref{PY2}(b) and Remark \ref{pairing-remark}(b) 
$v\Lambda_2(-M_2^{-1})\Lambda_2^T=\langle Q_+,Q\rangle_{F_2}=\eta(Q,K_2)$.
Therefore $C_L(v)=\beta+\eta(Q,K_2)$.
By Cochran's theorem (see Theorem \ref{cochran expansion})
$\eta(Q,K_2)=\sum_{n=1}^\infty(-1)^n\beta^n(Q,K_2)v^{2n}$.
By Lemma \ref{satolevine}(b) $(Q,K_2)$ is Cochran's derivative of $(K_1,K_2)$, so
$\beta^n(Q,K_2)=\beta^{n+1}(K_1,K_2)$. 
Also by Lemma \ref{satolevine}(a) $\beta=\beta^1(K_1,K_2)$.
Hence $C_L(v)=\beta+\eta(Q,K_2)=\sum_{n=0}^\infty(-1)^n\beta^{n+1}(K_1,K_2)v^{2n}$.
\end{proof}

In conclusion let us note a geometric criterion for the vanishing of Cochran's derived invariants
in a special case.

\begin{proposition} The following are equivalent for a PL link $(Q,K)$ with $\lk=0$:
\begin{enumerate}
\item all $\beta_i(Q,K)=0$ and $Q$ represents a conjugacy class in $\pi_1(S^3\but K)''$;
\item $Q$ bounds an immersed surface $\phi\:\Sigma\to S^3\but K$ such that $\phi^{-1}(Q)=Q$ 
and the image of every cycle in $\Sigma$ is null-homologous in $S^3\but K$.
\end{enumerate}
\end{proposition}

\begin{proof}
Let $M=S^3\but K$.
We have $\pi_1(M)''\simeq \pi_1(\tilde M)'$, where $\tilde M$ is the infinite cyclic cover of $M$
and its basepoint is some lift of the basepoint of $M$.
Let $\tilde Q$ and $\tilde Q_+$ be the nearby lifts in $\tilde M$ of $Q$ and of its parallel pushoff.

{\it (1)$\Rightarrow$(2).}
Since $Q$ represents a conjugacy class in $\pi_1(M)''$, $\tilde Q$ represents a conjugacy class in 
$\pi_1(\tilde M)'$, and therefore is null-homologous in $\tilde M$.
Then $\tilde Q$ bounds a PL embedded surface $\Sigma$ in $\tilde M$.
Since all $\beta_i(Q,K)=0$, we have $\eta_{(Q,K)}(t)=0$ and hence $\Sigma\cdot t^n\tilde Q_+=0$ for each $n$.
By attaching tubes to $\Sigma$ we may assume that $\Sigma\cap t^n\tilde Q_+=\emptyset$ for each $n$.
Then the restriction of the covering map $\tilde M\to M$ to $\Sigma$ is the desired immersion $\phi$.

{\it (2)$\Rightarrow$(1).}
Since the image of every cycle in $\Sigma$ is null-homologous in $S^3\but K$, $\phi$ lifts
to an immersion $\tilde\phi\:\Sigma\to\tilde M$.
It follows that $Q$ represents a conjugacy class in $\pi_1(S^3\but K)''$.
Since $\phi^{-1}(Q)=Q$, $\tilde\phi^{-1}(\bigcup_{n=\infty}^\infty t^nQ)=Q$.
Then also $\tilde\phi^{-1}(t^nQ_+)=\emptyset$ for each $n\ne 0$, and hence
$\eta_{(Q,K)}(t)$ contains only the constant term.
But $\eta_L(0)=1$ for every link $L$, so $\eta_{(Q,K)}(t)=0$.
Hence all $\beta_i(Q,K)=0$.
\end{proof}

\section*{Acknowledgements}

I would like to thank P. M. Akhmetiev, M. Il'insky, L. Traldi and A. Zastrow 
for useful discussions and thoughtful remarks. 
Additionally, the present work benefited from some of the student talks by N. Artyomov, M. Cabria, D. Gubarevich, 
M. Tyomkin and D. Zaytsev, who presented upon my suggestion some papers on the geometry of the Alexander 
polynomial at the Moscow Geometric Topology Seminar.

\end{document}